%% file: bmcomper.tex
\documentclass[11pt]{article}

\usepackage{amsmath}
\usepackage{amssymb}
\usepackage{color}
\usepackage{graphicx}
\usepackage{enumitem}
\usepackage{hyperref}
\usepackage{multirow}

\newcommand{\N}{\mathbb{N}}
\newcommand{\R}{\mathbb{R}}
\newcommand{\chico}{\mathrm{small}}
\newcommand{\ext}{\mathrm{ext}}
\newcommand{\ini}{\mathrm{ini}}

\newcommand{\new}{\mathrm{new}}

\newcommand{\refe}{\mathrm{ref}}
\newcommand{\SPG}{\mathrm{SPG}}
\newcommand{\target}{\mathrm{target}}
\newcommand{\temp}{\mathrm{temp}}
\newcommand{\trial}{\mathrm{trial}}
\newcommand{\opt}{\mathrm{opt}}
\newcommand{\feas}{\mathrm{feas}}
\newcommand{\compl}{\mathrm{compl}}

\newcommand{\best}{\mathrm{best}}
\newcommand{\halmos}{\hfill$\Box$}      
\newcommand{\half}{\frac{1}{2}}
\newcommand{\cbig}{c_{\mathrm{big}}}
\newcommand{\clips}{c_{\mathrm{lips}}}
\newcommand{\cinner}{c_{\mathrm{inner}}}
\newcommand{\cf}{c_f}
\newcommand{\deltalow}{\delta_{\mathrm{low}}}
\newcommand{\kend}{k_{\mathrm{end}}}
\newcommand{\rhobig}{\rho_{\mathrm{big}}}
\newcommand{\tol}{\mathrm{tol}}

\newtheorem{teo}{Theorem}[section]
\newtheorem{lem}{Lemma}[section]
\newenvironment{pro}{\noindent\textit{Proof:}}{\halmos\\}

\title{Complexity and performance of an \\ Augmented Lagrangian
  algorithm\footnote{This work was supported by FAPESP (grants
    2013/07375-0, 2016/01860-1, and 2018/24293-0) and CNPq (grants
    309517/2014-1 and 303750/2014-6).}}

\author{
  E. G. Birgin\thanks{Dept. of Computer Science, Institute of
    Mathematics and Statistics, University of S\~ao Paulo, Rua do
    Mat\~ao, 1010, Cidade Universit\'aria, 05508-090, S\~ao Paulo, SP,
    Brazil. email: egbirgin@ime.usp.br}
  \and
  J. M. Mart\'{\i}nez\thanks{Dept. of Applied Mathematics,
    Institute of Mathematics, Statistics, and Scientific Computing,
    State University of Campinas, 13083-859, Campinas, SP,
    Brazil. email: martinez@ime.unicamp.br}
}

\begin{document}

\date{July 4, 2019}

\maketitle

\begin{abstract}
Algencan is a well established safeguarded Augmented Lagrangian
algorithm introduced in [R. Andreani, E. G. Birgin,
  J. M. Mart\'{\i}nez and M. L. Schuverdt, On Augmented Lagrangian
  methods with general lower-level constraints, \textit{SIAM Journal
    on Optimization} 18, pp. 1286-1309, 2008]. Complexity results that
report its worst-case behavior in terms of iterations and evaluations
of functions and derivatives that are necessary to obtain suitable
stopping criteria are presented in this work. In addition, the
computational performance of a new version of the method is presented,
which shows that the updated software is a useful tool for solving
large-scale constrained optimization problems.\\

\end{abstract}

\section{Introduction}

Augmented Lagrangian methods have a long tradition in numerical
optimization.  The main ideas were introduced by Powell~\cite{powell},
Hestenes~\cite{hestenes}, and Rockafellar~\cite{rocka}. At each
(outer) iteration of an Augmented Lagrangian method one minimizes the
objective function plus a term that penalizes the non-fulfillment of
the constraints with respect to suitable shifted tolerances.  Whereas
the classical external penalty method~\cite{fiaccomccormick,fletcher}
needs to employ penalty parameters that tend to infinity, the shifting
technique aims to produce convergence by means of displacements of the
constraints that generate approximations to a solution with moderate
penalty parameters~\cite{bmbook}. As a by-product, one obtains
approximations of the Lagrange multipliers associated with the
original optimization problem. The safeguarded version of the
method~\cite{abmstango} discards Lagrange multiplier approximations
when they become very large. The convergence theory for safeguarded
Augmented Lagrangian methods was given
in~\cite{abmstango,bmbook}. Recently, examples that illustrate the
convenience of safeguarded Augmented Lagrangians were given
in~\cite{kanzowsteck}.

Conn, Gould, and Toint~\cite{cgtlancelot} produced the celebrated
package Lancelot, that solves constrained optimization problems using
Augmented Lagrangians in which the constraints are defined by
equalities and bounds. The technique was extended to the case of
equality constraints plus linear constraints
in~\cite{cgst}. Differently from Lancelot, in
Algencan~\cite{abmstango,bmbook} (see, also,
\cite{abmsobso,abmssecond,bfem,bfm,bmalgotan,bmfast,bmdecrease}), the
Augmented Lagrangian is defined not only with respect to equality
constraints but also with respect to inequalities. The theory
presented in~\cite{abmstango} and~\cite{bmbook} admits the presence of
lower-level constraints not restricted to boxes or polytopes. However,
in the practical implementations of Algencan, lower-level constraints
are always boxes.

In the last 10 years, the interest in Augmented Lagrangian methods was
renewed due to their ability to solve large-scale problems. Dost\'al
and Beremlijski~\cite{dostal,dostal2017} employed Augmented Lagrangian
methods for solving quadratic programming problems that appear in
structural optimization. Fletcher~\cite{fletcher2017} applied
Augmented Lagrangian ideas to the minimization of quadratics with box
constraints. Armand and Omheni~\cite{armand1} employed an Augmented
Lagrangian technique for solving equality constrained optimization
problems and handled inequality constraints by means of logarithmic
barriers~\cite{armand2}. Curtis, Gould, Jiang, and
Robinson~\cite{curtis1,curtis2} defined an Augmented Lagrangian
algorithm in which decreasing the penalty parameters is possible
following intrinsic algorithmic criteria. Local convergence results
without constraint qualifications were proved
in~\cite{fernandezsolodov}. The case with (possibly complementarity)
degenerate constraints was analyzed in \cite{ims}. Chatzipanagiotis
and Zavlanos~\cite{chatzi} defined and analyzed Augmented Lagrangian
methods in the context of distributed computation. An Exact Penalty
algorithm for constrained optimization with complexity results was
introduced in~\cite{cagt2011}. Grapiglia and Yuan \cite{gy2019}
analyzed the complexity of an Augmented Lagrangian algorithm for
inequality constraints based on the approach of Sun and Yuan
\cite{sunyuan} and assuming that a feasible initial point is
available.

In this paper, we report the main features of a new implementation of
Algencan.  The new Algencan preserves the main characteristics of the
previous algorithm: constraints are considered in the form of
equalities and inequalities, without slack variables; box-constrained
subproblems are solved using active-set strategies; and global
convergence properties are fully preserved. A new acceleration
procedure is introduced by means of which an approximate KKT point may
be obtained. It consists in applying a local Newton method to a
semismooth KKT system~\cite{mqi,qisun} starting from every Augmented
Lagrangian iterate. Special attention is given to the box-constraint
algorithm used for solving subproblems. The algorithm presented in
this paper is able to handle large-scale problems but not ``huge''
ones. This means that we deal with number of variables and Hessian
structures that make it affordable to use sparse
factorizations. Larger problems need the help of iterative linear
solvers which are not available in the new Algencan yet. Exhaustive
numerical experimentation is given and all the software employed is
available on a free basis in \url{http://www.ime.usp.br/~egbirgin/},
so that computational results are fully reproducible.

The paper is organized as follows. In Section~\ref{al}, we recall the
definition of Algencan with box lower-level constraints and we review
global convergence results. In Section~\ref{complexity}, we prove
complexity properties. In Section~\ref{newtonls}, we describe the
algorithm for solving box-constrained subproblems. In
Section~\ref{secimpl}, we describe the computer implementation. In
Section~\ref{secnumexp}, we report numerical experiments. Conclusions
are given in Section~\ref{secconcl}.\\

\noindent
\textbf{Notation.} If $C \subseteq \R^n$ is a convex set, $P_C(v)$
denotes the Euclidean projection of~$v$ onto~$C$. If $\ell, u \in
\R^n$, $[\ell,u]$ denotes the box defined by $\{ x \in \R^n \;|\; \ell
\leq x \leq u \}$. $(\cdot)_+ = \max\{0, \cdot \}$. If $v \in \R^n$,
$v_+$ denotes the vector with components $(v_i)_+$ for
$i=1,\dots,n$. If $v, w \in \R^n$, $\min\{v, w\}$ denotes the vector
with components $\min\{v_i, w_i\}$ for $i=1,\dots,n$. The symbol $\|
\cdot \|$ denotes the Euclidean norm.

\section{Augmented Lagrangian} \label{al}

In this section, we consider constrained optimization problems defined
by
\begin{equation} \label{nlp}
\mbox{Minimize } f(x) \mbox{ subject to } h(x) = 0, \; g(x) \leq 0, \mbox{ and } \ell \leq x \leq u,
\end{equation}
where $f: \R^n \to \R$, $h: \R^n \to \R^m$, and $g: \R^n \to \R^p$ are
continuously differentiable.

We consider the Augmented Lagrangian method in the way analyzed
in~\cite{abmstango} and \cite{bmbook}. This method has interesting
global theoretical properties. On the one hand, every limit point is a
stationary point of the problem of minimizing infeasibility. On the
other hand, every feasible limit point satisfies a sequential
optimality condition \cite{ahm,amrs1,amrs2}. This implies that every
feasible limit point is KKT-stationary under very mild constraint
qualifications \cite{amrs1,amrs2}. The basic definition of the method
and the main theoretical results are reviewed in this section.

The Augmented Lagrangian function~\cite{hestenes,powell,rocka}
associated with problem~(\ref{nlp}) is defined by
\[
L_{\rho}(x,\lambda,\mu) = f(x) +
\frac{\rho}{2} \left[ \sum_{i=1}^m \left( h_i(x) + \frac{\lambda_i}{\rho} \right)^2 +
  \sum_{i=1}^p \left( g_i(x) + \frac{\mu_i}{\rho} \right)_+^2 \right]
\]
for all $x \in [\ell,u]$, $\rho > 0$, $\lambda \in \R^m$, and $\mu \in
\R^p_+$. The Augmented Lagrangian model algorithm
follows.\\

\noindent
\textbf{Algorithm~\ref{al}.1:} Assume that $x^0 \in \R^n$,
$\lambda_{\min} < \lambda_{\max}$, ${\bar \lambda}^1 \in
       [\lambda_{\min}, \lambda_{\max}]^m$, $\mu_{\max} > 0$, ${\bar
         \mu}^1 \in [0, \mu_{\max}]^p$, $\rho_1 > 0$, $\gamma > 1$, $0
       < \tau < 1$, and $\{\varepsilon_k\}_{k=1}^\infty$ are
       given. Initialize $k \leftarrow 1$.

\begin{description}

\item[Step 1.] Find $x^k \in [\ell,u]$ as an approximate solution to
  \begin{equation} \label{subprob}
  \mbox{Minimize } L_{\rho_k}(x,{\bar \lambda}^k, {\bar \mu}^k) 
  \mbox{ subject to } \ell \leq x \leq u
  \end{equation}
  satisfying
  \begin{equation} \label{subprostop}
  \left\| P_{[\ell,u]}\left(x^k - \nabla L_{\rho_k}(x^k,{\bar \lambda}^k,{\bar \mu}^k) \right) - x^k \right\|
  \leq \varepsilon_k.
  \end{equation}
  
\item[Step 2.] Define
  \[
  V^k = \min \left\{ -g(x^k), \frac{{\bar \mu}^k}{\rho_k} \right\}.
  \]
  If $k = 1$ or
  \begin{equation} \label{testfeas} 
  \max \left\{ \|h(x^k)\|, \|V^k\| \right\} \leq
  \tau \max \left\{ \|h(x^{k-1})\|, \|V^{k-1}\| \right\},
  \end{equation}
  choose $\rho_{k+1} = \rho_k$. Otherwise, define $\rho_{k+1} = \gamma
  \rho_k$.

\item[Step 3.] Compute
  \begin{equation} \label{lambdamas} 
  \lambda^{k+1} = {\bar \lambda}^k + \rho_k h(x^k)
  \mbox{ and }
  \mu^{k+1} = \left( {\bar \mu}^k + \rho_k g(x^k) \right)_+.
  \end{equation}
  Compute ${\bar \lambda}^{k+1} \in [\lambda_{\min},
    \lambda_{\max}]^m$ and ${\bar \mu}^{k+1}_i \in [0, \mu_{\max}]^p$.
  Set $k \leftarrow k+1$ and go to Step~1.

\end{description}

The problem of finding an approximate minimizer of $L_{\rho_k}(x,{\bar
  \lambda}^k, {\bar \mu}^k)$ onto $[\ell,u]$ in the sense
of~(\ref{subprostop}) can always be solved. In fact, due to the
compactness of $[\ell,u]$, a global minimizer, that obviously
satisfies~(\ref{subprostop}), always exists. Moreover, local
minimization algorithms are able to find an approximate stationary
point satisfying~(\ref{subprostop}) in a finite number of
iterations. Therefore, given an iterate~$x^k$, the iterate~$x^{k+1}$
is well defined. So, Algorithm~\ref{al}.1 generates an infinite
sequence $\{x^k\}$ whose properties are surveyed below. Of course, as
it will be seen later, suitable stopping criteria can be defined by
means of which acceptable approximate solutions to~(\ref{nlp}) are
usually obtained.

Algorithm~\ref{al}.1 has been presented without a ``stopping
criterion''.  This means that, in principle, the algorithm generates
an infinite sequence of primal iterates $x^k$ and Lagrange-multiplier
estimators. Complexity results presented in this work report the
worst-case effort that could be necessary to obtain different
properties, that may be used as stopping criteria in practical
implementations or not. We believe that the interpretation of these
results helps to decide which stopping criteria should be used in a
practical application.

The relevant theoretical properties of this algorithm are the
following:
\begin{enumerate}
\item Every limit point $x^*$ of the sequence generated by the
  algorithm satisfies the complementarity condition
  \begin{equation} \label{complementarity}
    \mu^{k+1}_i = 0 \mbox{ whenever } g_i(x^*) < 0
  \end{equation}
  for $k$ large enough. (See~\cite[Thm.4.1]{bmbook}.)

\item Every limit point $x^*$ of the sequence generated by the
  algorithm satisfies the first-order optimality conditions of the
  feasibility problem
  \begin{equation} \label{fisipro}
    \mbox{Minimize } \|h(x)\|^2 + \|g(x)_+\|^2 \mbox{ subject to } \ell \leq x \leq u.
  \end{equation}
  (See~\cite[Thm.6.5]{bmbook}.)

\item If, for all $k \in \{1, 2, \dots\}$, $x^k$ is an approximate
  global minimizer of $L_{\rho_k}(x, \bar \lambda^k, \bar \mu^k)$ onto
  $[\ell,u]$ with tolerance $\eta > 0$, every limit point of $\{x^k\}$
  is a global minimizer of the infeasibility function $\|h(x)\|^2 +
  \|g(x)_+\|^2$. Condition~(\ref{subprostop}) does not need to hold in
  this case. (See~\cite[Thm.5.1]{bmbook}.)

\item If, for all $k \in \{1, 2, \dots \}$, $x^k$ is an approximate
  global minimizer of $L_{\rho_k}(x,\bar \lambda^k,\bar \mu^k)$ onto
  $[\ell,u]$ with tolerance $\eta_k \downarrow 0$, every feasible
  limit point of $\{x^k\}$ is a global minimizer of the general
  constrained minimization problem~(\ref{nlp}). As before,
  condition~(\ref{subprostop}) is not necessary in this
  case. (See~\cite[Thm.5.2]{bmbook}.)

\item If $\varepsilon_k \downarrow 0$, every feasible limit point of
  the sequence $\{x^k\}$ satisfies the sequential optimality
  condition~AKKT \cite{ahm} given by
  \begin{equation} \label{akktproj}
    \lim_{k \in K} \left\| P_{[\ell,u]} \left( x^k - \left( \nabla
      f(x^k) + \nabla h(x^k) \lambda^{k+1} + \nabla g(x^k) \mu^{k+1}
      \right) \right) -x^k \right\| = 0
  \end{equation}
  and
  \begin{equation} \label{compleme}
    \lim_{k \in K} \max \{ \|h(x^k)\|_\infty,  \| \min\{-g(x^k), \mu^{k+1}\}\|_\infty \} = 0,
  \end{equation}
  where the sequence of indices $K$ is such that $\lim_{k \in K} x^k =
  x^*$. (See~\cite[Thm.6.4]{bmbook}.)
\end{enumerate}
Under an additional Lojasiewicz-like condition, it is obtained that
$\lim_{k \in K} \sum_{i=1}^p \mu^{k+1}_i g_i(x^k) = 0$
(see~\cite{amscakkt}).  Moreover, in \cite{afss}, it was proved that
an even stronger sequential optimality condition is satisfied by the
sequence $\{x^k\}$, perhaps associated with different Lagrange
multipliers approximations than the ones generated by the Augmented
Lagrangian algorithm.

These properties say that, even if $\varepsilon_k$ does not tend to
zero, Algorithm~\ref{al}.1 finds stationary points of the
infeasibility measure $\|h(x)\|^2 + \|g(x)_+\|^2$ and that, when
$\varepsilon_k$ tends to zero, feasible limit points satisfy a
sequential optimality condition. Thus, under very weak constraint
qualifications, feasible limit points satisfy Karush-Kuhn-Tucker
conditions. See \cite{amrs1,amrs2}. Some of these properties, but not
all, are shared by other constrained optimization algorithms. For
example, the property that feasible limit points satisfy optimality
KKT conditions is proved to be satisfied by other optimization
algorithms only under much stronger constraint qualifications than the
ones required by Algorithm~\ref{al}.1. Moreover, the Newton-Lagrange
method may fail to satisfy approximate KKT conditions even when it
converges to the solution of rather simple constrained optimization
problems \cite{amsfail,amss}.

Augmented Lagrangian implementations have a modular structure. At each
iteration, a box-constrained optimization problem is approximately
solved. The efficiency of the Augmented Lagrangian algorithm is
strongly linked to the efficiency of the box-constraint solver.

Algencan may be considered to be a conservative visit to the Augmented
Lagrangian framework. For example, subproblems are solved with
relatively high precision, instead of stopping subproblem solvers
prematurely according to information related to the constrained
optimization landscape. It could be argued that solving subproblems
with high precision at points that may be far from the solution
represents a waste of time. Nevertheless, our point of view is that
saving subproblem iterations when one is close to a subproblem
solution is not worthwhile because in that region Newton-like solvers
tend to be very fast; and accurate subproblems' solutions help to
produce better approximations of Lagrange multipliers. Algencan is
also conservative when subproblems' solvers use minimal information
about the structure of the Augmented Lagrangian function they
minimize. The reason for this decision is connected to the modular
structure of Algencan. Subproblem solvers are continuously being
improved due to the continuous and fruitful activity in
bound-constrained minimization. Therefore, we aim to take advantage of
those improvements with minimal modifications of subproblem algorithms
when applied to minimize Augmented Lagrangians.

\section{Complexity} \label{complexity}

This section is devoted to worst-case complexity results related to
Algorithm~\ref{al}.1. Algorithm~\ref{al}.1 was not devised with the
aim of optimizing complexity. Nevertheless, our point of view is that
the complexity analysis that follows helps to understand the actual
behavior of the algorithm, filling a gap opened by the convergence
theory.
  
By~(\ref{lambdamas}) and straightforward calculations, we have that,
for all $k = 1, 2, 3, \dots$,
\[
\nabla f(x^k) + \nabla h(x^k) \lambda^{k+1} + \nabla g(x^k) \mu^{k+1}
= \nabla L_{\rho_k}(x^k, \bar \lambda^k, \bar \mu^k).
\]
Therefore, the fulfillment of
\begin{equation} \label{paraparar1}
\|P_{[\ell, u]} (x^k - \nabla  L_{\rho_k}(x^k, \bar \lambda^k, \bar \mu^k)) - x^k) \|
\leq \varepsilon
\end{equation}
implies that the projected gradient of the Lagrangian with
multipliers~$\lambda^{k+1}$ and~$\mu^{k+1}$ approximately vanishes
with precision~$\varepsilon$. In the next lemma, we show that the
fulfillment of
\begin{equation} \label{paraparar2}
\max \{ \|h(x^k)\|_\infty, \|V_k\|_\infty\} \leq \delta
\end{equation} 
implies that feasibility and complementarity hold at~$x^k$ with
precision~$\delta$. For these reasons, in the context of
Algorithm~\ref{al}.1, iterates that satisfy~(\ref{paraparar1})
and~(\ref{paraparar2}) are considered approximate stationary points of
problem~(\ref{nlp}).

\begin{lem} \label{lemcomplexity1}
For all $\delta > 0$,   
\begin{equation} \label{tesfeacom}
\max \{\|h(x^k)\|_\infty, \|V_k\|_\infty\} \leq \delta
\end{equation}  
implies that
\begin{equation} \label{lastres}
  \|h(x^k)\|_\infty \leq \delta, \;
  \|g(x^k)_+\|_\infty \leq \delta, \mbox{ and},
  \mbox{ for all } j=1,\dots,p, \mu^{k+1}_j = 0
  \mbox{ if } g_j(x^k) < - \delta.
\end{equation}
\end{lem}

\begin{pro}
By~(\ref{tesfeacom}), $\|h(x^k)\|_\infty \leq \delta$ and
$|\min\{-g_j(x^k), \bar{\mu}_j^k/\rho_k\}| \leq \delta$ for all
$j=1,\dots,p$. Therefore, $-g_j(x^k) \geq -\delta$, so $g_j(x^k) \leq
\delta$ for all $j=1\dots, p$. Moreover, by (\ref{tesfeacom}), if
$g_j(x^k) < -\delta$, we necessarily have that $\bar{\mu}_j^k/\rho_k
\leq \delta$. Adding these two inequalities, we obtain that, if
$g_j(x^k) < -\delta$ then $g_j(x^k) + \bar{\mu}_j^k /\rho_k <
0$. Consequently, $\rho_k g_j(x^k) + \bar{\mu}_j^k < 0$, so
$\mu^{k+1}_j = 0$. Therefore, (\ref{tesfeacom})
implies~(\ref{lastres}) as we wanted to prove.
\end{pro}

In Theorem~\ref{complexity}.1 below, we assume that the sequence
$\{\rho_k\}$ is bounded. Sufficient conditions for this requirement,
where the bound $\bar{\rho}$ only depends on algorithmic parameters
and characteristics of the problem, were given in \cite{abmstango} and
\cite{bmbook}. We also assume that there exists $N(\varepsilon) \in
\{1, 2, 3, \dots\}$ such that $\varepsilon_k \leq \varepsilon$ for all
$k \geq N(\varepsilon)$. Clearly, this condition can be enforced by
the criterion used to define~$\{\varepsilon_k\}$. For example,
$\varepsilon_{k+1} = \half \varepsilon_k$ obviously implies that
$\varepsilon_k \leq \varepsilon$ if $k > N(\varepsilon) \equiv
\log(\varepsilon)/\log(\varepsilon_1)$.

\begin{lem} \label{lemcomplexity2}
There exists $\cbig > 0$ such that, for all $k \geq 1$,
\begin{equation} \label{defcbig}
\max \{\|h(x^k)\|_\infty, \|V_k\|_\infty\} \leq \cbig.
\end{equation} 
\end{lem}

\begin{pro}
Since, by definition of the algorithm, $\rho_k \geq \rho_1$, the
bound~(\ref{defcbig}) comes from the continuity of $h$ and $g$, the
compactness of the domain $[\ell, u]$, and the boundedness of $\bar
\mu^k$.
\end{pro}

From now on, $\cbig$ will denote a positive constant
satisfying~(\ref{defcbig}), whose existence is guaranteed by
Lemma~\ref{lemcomplexity2}.

\begin{teo} \label{teocomplexity1}
Let $\delta > 0$ and $\varepsilon > 0$ be given. Assume that, for all
$k \in \{1, 2, 3, \dots\}$, $\rho_k \leq \bar{\rho}$. Moreover, assume
that, for all $k \geq N(\varepsilon)$, we have that $\varepsilon_k
\leq \varepsilon$. Then, after at most
\begin{equation} \label{iteraciones}
N(\varepsilon) + \left[ \log(\bar{\rho}/\rho_1)/\log(\gamma)
  \right] \times \left[ \log(\delta/\cbig)/\log (\tau) \right]
\end{equation}
iterations, we obtain $x^k \in [\ell, u]$, $\lambda^{k+1} \in \R^m$,
and $\mu^{k+1} \in \R^p_+$ such that
\begin{equation} \label{aprokkt1}
\left\| P_{[\ell,u]} \left( x^k - \left( \nabla f(x^k) + \nabla h(x^k)
  \lambda^{k+1} + \nabla g(x^k) \mu^{k+1} \right) \right) -x^k
\right\| \leq \varepsilon,
\end{equation}
\begin{equation} \label{aprokkt2}
  \| h(x^k)\|_\infty \leq \delta, \;
  \|g(x^k)_+\|_\infty \leq \delta, 
\end{equation}
and, for all $j=1,\dots,p$, 
\begin{equation} \label{aprokkt3}
\mu^{k+1}_j = 0 \mbox{ whenever } g_j(x^k) <  - \delta.
\end{equation}                                                                                                   
\end{teo}

\begin{pro}
The number of iterations such that $\rho_{k+1} =
\gamma \rho_k$ is bounded above by
\begin{equation} \label{bounrho}  
\log (\bar{\rho}/ \rho_1)/\log(\gamma).
\end{equation}
Therefore, this is also a bound for the number of iterations at
which~(\ref{testfeas}) does not hold.

By (\ref{defcbig}), if (\ref{testfeas}) holds during
\begin{equation}\label{bounconse}
\log (\delta/\cbig)/\log \tau
\end{equation}
consecutive iterations, we get that   
\[
\max \{\|h(x^k)\|_\infty, \|V_k\|_\infty\} \leq \delta,
\]
which, by Lemma~\ref{lemcomplexity1}, implies (\ref{aprokkt2})
and~(\ref{aprokkt3}).

Now, by hypothesis, after $N(\varepsilon)$ iterations, we have that
$\varepsilon_k \leq \varepsilon$. Therefore, by~(\ref{bounrho})
and~(\ref{bounconse}), after at most
\begin{equation} \label{iteracionesN}
N(\varepsilon) + [ \log (\bar{\rho}/ \rho_1)/\log(\gamma)] \times
[\log (\delta/\cbig)/\log (\tau)]
\end{equation}
iterations, we have that~(\ref{aprokkt1}), (\ref{aprokkt2}),
and~(\ref{aprokkt3}) hold.
\end{pro}

Theorem~\ref{teocomplexity1} shows that, as expected, if $\rho_k$ is
bounded, we obtain approximate feasibility and optimality. In the
following theorem, we assume that the subproblems are solved by means
of some method that, for obtaining precision $\varepsilon > 0$,
employs at most $c \varepsilon^{-q}$ iterations and evaluations, where
$c$ only depends on characteristics of the problem, the upper bound
for $\rho_k$, and algorithmic parameters of the method.

\begin{teo} \label{teocomplexity2}
In addition to the hypotheses of Theorem~\ref{teocomplexity1}, assume
that there exist $\cinner > 0$ and $q > 0$, where $\cinner$ only
depends on~$\bar{\rho}$, $\lambda_{\min}$, $\lambda_{\max}$,
$\mu_{\max}$, $\ell$, $u$, and characteristics of the functions~$f$,
$h$, and~$g$, such that the number of inner iterations, function and
derivative evaluations that are necessary to obtain (\ref{subprostop})
is bounded above by $\cinner \varepsilon_k^{-q}$. Then, the number of
inner iterations, function evaluations, and derivative evaluations
that are necessary to obtain~$k$ such that~(\ref{aprokkt1}),
(\ref{aprokkt2}), and~(\ref{aprokkt3}) hold is bounded above by
\[
\cinner \varepsilon_{\min}^{-q} \left\{ N(\varepsilon) + \left[ \log(\bar{\rho}/\log
  \rho_1)/\log(\gamma) \right] \times \left[\log (\delta/\cbig)/\log
  (\tau) \right] \right\},
\]    
where
\begin{equation} \label{epsilonmin}
\varepsilon_{\min} = \min\{\varepsilon_k \;|\; k \leq N(\varepsilon) + \left[ \log(\bar{\rho}/\log
  \rho_1)/\log(\gamma) \right] \times \left[\log (\delta/\cbig)/\log
  (\tau) \right] \}.
\end{equation}     
\end{teo}

\begin{pro}
The desired result follows from Theorem~\ref{complexity}.1 and the
assumptions of this theorem.
\end{pro}

Note that, in Theorem~\ref{teocomplexity2}, we admit the possibility
that $\varepsilon_k$ decrease after completing $N(\varepsilon)$
iterations. This is the reason for the definition of
$\varepsilon_{\min}$ (\ref{epsilonmin}). In practical implementations,
it is reasonable to stop decreasing $\varepsilon_k$ when it achieves a
user-given stopping tolerance $\varepsilon$.  According to
Theorem~\ref{teocomplexity2}, the complexity bounds related to
approximate optimality, feasibility, and complementarity depend on the
optimality tolerance $\varepsilon$ in, essentially, the same way that
the complexity of the subproblem solver depends on its stopping
tolerance. In other words, under the assumption of boundedness of
penalty parameters, the worst-case complexity of the Augmented
Lagrangian method is essentially the same as the complexity of the
subproblem solver.

In computer implementations, it is usual to employ, in addition to a
(successful) stopping criterion based on (\ref{aprokkt1}),
(\ref{aprokkt2}), and (\ref{aprokkt3}), an (unsuccessful) stopping
criterion based on the size of the penalty parameter.  The rationale
is that if the penalty parameter grew to be very large, it is not
worthwhile to expect further improvements with respect to feasibility
and we are probably close to an infeasible local minimizer of
infeasibility. The complexity results that correspond to this decision
are given below.

\begin{teo} \label{teocomplexitygrande}
Let $\delta > 0$, $\varepsilon > 0$, and $\rhobig > \rho_1 $ be
given. Assume that, for all $k \geq N(\varepsilon)$, we have that
$\varepsilon_k \leq \varepsilon$. Then, after at most
\begin{equation} \label{iteracoes}
N(\varepsilon) + \left[ \log(\rhobig/\rho_1)/\log(\gamma)
  \right] \times \left[ \log(\delta/\cbig)/\log (\tau) \right]
\end{equation}
iterations, we obtain $x^k \in [\ell, u]$, $\lambda^{k+1} \in \R^m$,
and $\mu^{k+1} \in \R^p_+$ such that (\ref{aprokkt1}),
(\ref{aprokkt2}), and (\ref{aprokkt3}) hold or we obtain an iteration
such that $\rho_k \geq \rhobig$.                                                               
\end{teo}

\begin{pro}
If $\rho_k \leq \rhobig$ for all $k \leq N(\varepsilon) + \left[
  \log(\rhobig/\rho_1)/\log(\gamma) \right] \times \left[
  \log(\delta/\cbig)/\log (\tau) \right]$, by the same argument used
in the proof of Theorem~\ref{teocomplexity1}, with $\rhobig$
replacing $\bar{\rho}$, we obtain that (\ref{aprokkt1}),
(\ref{aprokkt2}), and (\ref{aprokkt3}) hold.
\end{pro}

\begin{teo} \label{teocomplexitygrande2}
In addition to the hypotheses of Theorem~\ref{teocomplexitygrande},
assume that there exist $\cinner > 0$ and $q > 0$, where $\cinner$
only depends on~$\rhobig$, $\lambda_{\min}$, $\lambda_{\max}$,
$\mu_{\max}$, $\ell$, $u$, and characteristics of the functions~$f$,
$h$, and~$g$, such that the number of inner iterations, function and
derivative evaluations that are necessary to obtain (\ref{subprostop})
is bounded above by $\cinner \varepsilon_k^{-q}$. Then, the number of
inner iterations, function evaluations, and derivative evaluations
that are necessary to obtain~$k$ such that~(\ref{aprokkt1}),
(\ref{aprokkt2}), and~(\ref{aprokkt3}) hold or such that $\rho_k >
\rhobig$ is bounded above by
\[
\cinner \varepsilon_{\min,2}^{-q} \left\{ N(\varepsilon) + \left[ \log(\rhobig/
  \rho_1)/\log(\gamma) \right] \times \left[\log (\delta/\cbig)/\log
  (\tau) \right] \right\},
\]    
where
\begin{equation} \label{2epsilonmin}
\varepsilon_{\min,2} = \min\{\varepsilon_k \;|\; k \leq \left\{
N(\varepsilon) + \left[ \log(\rho_{\max}/\log \rho_1)/\log(\gamma)
  \right] \times \left[\log (\delta/\cbig)/\log (\tau) \right]
\right\}.
\end{equation}     
\end{teo}

\begin{pro}
The desired result follows directly from Theorem~\ref{teocomplexitygrande}.
\end{pro}

The complexity results proved up to now indicate that suitable
stopping criteria for Algorithm~\ref{al}.1 could be based on the
fulfillment of (\ref{aprokkt1}), (\ref{aprokkt2}), and
(\ref{aprokkt3}) or, alternatively, on the occurrence of an
undesirable big penalty parameter. The advantage of these criteria is
that, according to them, worst-case complexity is of the same order as
the complexity of subproblem solvers. Convergence results establish
that solutions obtained with very large penalty parameters are close
to stationary points of the infeasibility. However, stationary points
of infeasibility may be feasible points and, again, convergence theory
shows that when Algorithm~\ref{al}.1 converges to a feasible point,
this point satisfies AKKT optimality conditions, independently of
constraint qualifications. As a consequence, the danger exists of
interrupting executions prematurely, in situations in which meaningful
progress could be obtained admitting further increases of the penalty
parameter. This state of facts leads one to analyze complexity of
Algorithm~\ref{al}.1 independently of penalty parameter growth and
introducing a possibly more reliable criterion for detecting
infeasible stationary points of infeasibility. Roughly speaking, we
will say that an iterate seems to be an infeasible stationary point of
infeasibility when the projected gradient of the infeasibility measure
is significantly smaller than the infeasibility value. The natural
question that arises is whether the employment of this (more reliable)
stopping criterion has an important effect on the complexity bounds.

Assumptions on the limitation of $\rho_k$ are given up from now on.
Note that the possibility that $\rho_k \to \infty$ needs to be
considered since necessarily takes place, for example, when the
feasible region is empty.

\begin{lem} \label{lemcomplexity3}
There exist $\clips$, $\cf > 0$ such that, for all $x \in [\ell, u]$,
$\lambda \in [\lambda_{\min}, \lambda_{\max}]^m$, and $\mu \in [0,
  \mu_{\max}]^p$, one has
\begin{equation} \label{clips}
\|\nabla h(x)\| \|\lambda\| + \|\nabla g(x)\| \|\mu\|  \leq \clips
\end{equation}  
and     
\begin{equation} \label{cf}
\|\nabla f(x)\| \leq \cf.
\end{equation}
\end{lem}

\begin{pro}
The desired result follows from the boundedness of the domain, the
continuity of the functions, and the boundedness of~$\lambda$
and~$\mu$.
\end{pro}

The following lemma establishes a bound for the projected gradient of
the infeasibility measure in terms of the value of the displaced
infeasibility and the value of the penalty parameter.

\begin{lem} \label{lemcomplexity4}
For all $x \in [\ell, u]$, $\lambda \in [\lambda_{\min},
  \lambda_{\max}]^m$, $\mu \in [0, \mu_{\max}]^p$, and $\rho > 0$, one
has that
\[
\left\|P_{[\ell, u]}\left(x - \nabla\left[ \|h(x)\|^2] + \|g(x)_+\|^2 \right] \right) - x \right\|
\]
\[
\leq \left\|  P_{[\ell,u]}\left(x - \nabla\left[\|h(x) + \lambda/\rho\|^2 + 
\| (g(x)+\mu/\rho)_+ \|^2 \right] \right) - x \right\| +  2 \clips / \rho,
\]
where $\clips$ is defined in Lemma~\ref{lemcomplexity3}.
\end{lem}

\begin{pro}
Note that 
\[
\half \nabla \left[ \left\|h(x) + \lambda/\rho \right\|^2 + \left\| \left( g(x)+\mu/\rho \right)_+ \right\|^2 \right]
= h'(x)^T \left( h(x) + \lambda/\rho \right) + g'(x)^T \left( g(x) + \mu/\rho \right)_+
\]
and
\[
\half \nabla \left[ \|h(x)\|^2 + \|g(x)_+\|^2 \right] =
\nabla h(x) h(x) + \nabla g(x) g(x)_+.
\]   
Therefore,
\[
\left\|\half \nabla \left[ \| h(x) + \lambda/\rho \|^2 +
\| \left( g(x) + \mu/\rho)_+ \right\|^2 \right] -
\half \nabla \left[ \|h(x)\|^2 + \|g(x)_+\|^2 \right] \right\|
\]     
\[
\leq
\left\| \nabla h(x) \lambda/\rho + \nabla g(x)
\left[ (g(x) + \mu/\rho)_+ - g(x)_+ \right] \right\|
\leq
\frac{1}{\rho} \left[ \| \nabla h(x) \| \|\lambda\| +
\| \nabla g(x)\| \|\mu\| \right].
\]
Then, by~(\ref{clips}), if $\rho > 0$, $x \in [\ell, u]$, $\lambda \in
[\lambda_{\min}, \lambda_{\max}]^m$, and $\mu \in [0, \mu_{\max}]^p$,
\[
\left\| \nabla \left[ \|h(x)\|^2 + \|g(x)_+\|^2 \right] - \nabla
  \left[ \|h(x) + \lambda/\rho\|^2 + \| \left( g(x)+\mu/\rho
    \right)_+\|^2 \right] \right\| \leq 2 \clips / \rho.
\]  
So, by the non-expansivity of projections,
\[
\resizebox{\textwidth}{!}{$
\left\| P_{[\ell, u]}\left(x - \nabla \left[ \|h(x)\|^2 + \|g(x)_+\|^2 \right]\right) -
P_{[\ell,u]}\left(x - \nabla \left[ \|h(x) + \lambda/\rho\|^2 +
  \| (g(x)+\mu/\rho)_+ \|^2 \right] \right) \right\| \leq 2 \clips/\rho.
$}
\]
Thus, the thesis is proved.
\end{pro}

The following theorem establishes that, before the number of
iterations given by~(\ref{iterteo5}), we necessarily find an
approximate KKT point or we find an infeasible point that, very
likely, is close to an infeasible stationary point of the
infeasibility measure. The latter type of infeasible points is
characterized by the fact that the projected gradient of the
infeasibility is smaller than~$\deltalow$ whereas the infeasibility
value is bigger than $\delta \gg \deltalow$.
           
\begin{teo} \label{teocomplexity3}
Let $\delta > 0$, $\deltalow \in (0,\delta)$, and $\varepsilon > 0$ be
given. Assume that $N(\deltalow, \varepsilon)$ is such that
$\varepsilon_k \leq \min\{\varepsilon, \deltalow \}/4$ for all $k \geq
N(\deltalow, \varepsilon)$. Then, after at most
\begin{equation} \label{iterteo5}
N(\deltalow, \varepsilon) + \left[ \frac{\log
    (\delta/\cbig)}{\log(\tau)}\right] \times \left[ \frac{\log \left(
    \rho_{\max} / \rho_1 \right)}{\log(\gamma)} \right]
\end{equation}      
iterations, where
\begin{equation} \label{rhomax}
\rho_{\max} = \max\left\{1, \frac{4 \clips}{\deltalow},
\frac{\mu_{\max}}{\delta}, \frac{4 c_f}{\deltalow } \right\},
\end{equation}
we obtain an iteration $k$ such that one of the following
two facts takes place:
\begin{enumerate}
\item The iterate $x^k \in [\ell, u]$ verifies
  \begin{equation} \label{paradamala}
    \resizebox{0.9\textwidth}{!}{$
    \left\| P_{[\ell,u]}\left( x^k - \nabla \left[ \|h(x^k)\|^2 + \|g(x^k)_+\|^2 \right] \right) - x^k \right\|
    \leq \deltalow \mbox{ and } \max\{\|h(x^k)\|_\infty, \|g(x^k)_+\|_\infty\}> \delta.
    $}
  \end{equation}

\item The multipliers $\lambda^{k+1} \in \R^m$ and $\mu^{k+1} \in
  \R^p_+$ are such that
  \begin{equation} \label{primaldual1}
    \left\| P_{[\ell,u]}\left( x^k - \left( \nabla f(x^k) + \nabla h(x^k)
      \lambda^{k+1} + \nabla g(x^k) \mu^{k+1} \right) \right) -x^k \right\|
    \leq \varepsilon,
  \end{equation}
  \begin{equation} \label{primaldual2}
    \|h(x^k)\|_\infty \leq \delta, \; \|g(x^k)_+\|_\infty \leq \delta, 
  \end{equation}
  and, for all $j=1,\dots,p$, 
  \begin{equation} \label{primaldual3}
    \mu^{k+1}_j = 0 \mbox{ whenever } g_j(x^k) <  - \delta.
  \end{equation}
\end{enumerate}
\end{teo}

\begin{pro}
Let $\kend$ be such that
\begin{equation} \label{implica}
\left\| P_{[\ell,u]}\left( x^k - \nabla \left[ \|h(x^k)\|^2 + \|g(x^k)_+\|^2 \right] \right) - x^k \right\|
\leq \deltalow \Rightarrow \max \{ \|h(x^k)\|_\infty, \|g(x^k)_+\|_\infty \} \leq \delta
\end{equation}
for all $k \leq \kend$ whereas~(\ref{implica}) does not hold if $k =
\kend + 1$. (With some abuse of notation, we say that $\kend = \infty$
when~(\ref{implica}) holds for all $k$.) In other words, if $k \leq
\kend$,
\begin{equation} \label{disyun}
\left\| P_{[\ell,u]}\left( x^k - \nabla \left[ \|h(x^k)\|^2 + \|g(x^k)_+\|^2 \right] \right) - x^k \right\|
> \deltalow \mbox{ or } \max \{ \|h(x^k)\|_\infty, \|g(x^k)_+\|_\infty \} \leq \delta,
\end{equation}    
whereas~(\ref{disyun}) does not hold if $k = \kend + 1$.

We consider two possibilities:
\begin{equation} \label{kendmenor}
\kend <
N(\deltalow, \varepsilon) + \left[ \frac{\log
    (\delta/\cbig)}{\log(\tau})\right] \times \left[ \frac{\log \left(
    \rho_{\max} / \rho_1 \right)}{\log(\gamma)} \right]
\end{equation}  
and
\begin{equation} \label{kendmayor}
\kend \geq
N(\deltalow, \varepsilon) + \left[ \frac{\log
    (\delta/\cbig)}{\log(\tau})\right] \times \left[ \frac{\log \left(
    \rho_{\max} / \rho_1 \right)}{\log(\gamma)} \right].
\end{equation}      
In the first case, since~(\ref{implica}) does not hold for $k =
\kend+1$, it turns out that~(\ref{paradamala}) occurs at iteration
$\kend + 1$. It remains to analyze the case in which~(\ref{kendmayor})
takes place.

Suppose that
\begin{equation} \label{kbajo}
k \leq
N(\deltalow, \varepsilon) + \left[ \frac{\log
    (\delta/\cbig)}{\log(\tau})\right] \times \left[ \frac{\log \left(
    \rho_{\max} / \rho_1 \right)}{\log(\gamma)} \right],
\end{equation}
\begin{align}
\varepsilon_k &\leq \deltalow / 4, \label{coneps} \\[2mm]
\rho_k &\geq 1, \label{rhouno} \\[2mm]
\rho_k &\geq 4 c_f / \deltalow,  \label{rhocf} \\[2mm]
\rho_k &\geq 4 \clips / \deltalow, \label{rholis} \\[2mm]
\rho_k &\geq \mu_{\max} / \delta, \label{rhomuma} \\[2mm]
k &\geq N(\deltalow, \varepsilon). \label{kene}
\end{align}
By (\ref{subprostop}), for all $k \geq 1$, we have that
\[
\left\| P_{[\ell,u]}\left( x^k - \nabla f(x^k) - \frac{\rho_k}{2}
\nabla \left\{ \sum_{i=1}^m \left[ h_i(x^k) +
  \frac{\bar{\lambda}_i^k}{\rho_k}\right ]^2 + \sum_{i=1}^p \left[ \left(
  g_i(x^k) + \frac{\bar{\mu}_i^k}{\rho_k} \right)_+ \right]^2
\right\} \right) - x^k \right\| \leq \varepsilon_k.
\]   
Therefore, by (\ref{rhouno}),  
\[
\left\|
P_{[\ell,u]}\left(
x^k - \frac{1}{\rho_k} \nabla f(x^k) -
\half \nabla \left(
\| h(x^k) + \bar \lambda^k/\rho_k \|^2 +
\| (g(x^k) + \bar \mu^k/\rho_k)_+ \|^2
\right)
\right) - x^k
\right\| \leq \varepsilon_k.
\] 
Therefore, by the non-expansivity of projections and (\ref{cf}), we
have that
\begin{equation} \label{lafalta}
  \left\| P_{[\ell,u]}\left( x^k - \half \nabla \left( \|h(x^k) +
  \bar \lambda^k/\rho_k\|^2 + \|(g(x^k) + \bar \mu^k/\rho_k)_+\|^2 \right) \right) - x^k
  \right\| \leq \varepsilon_k + \frac{c_f}{\rho_k}.
\end{equation} 
So, by (\ref{coneps}) and (\ref{rhocf}), 
\begin{equation} \label{quiero1}
  \left\| P_{[\ell,u]}\left(x^k - \nabla\left( \|h(x^k) + \bar \lambda^k/\rho_k\|^2
    + \|(g(x^k)+\bar \mu^k/\rho_k)_+\|^2\right) \right) - x^k \right\| \leq \deltalow/2.
\end{equation}     
Therefore, by Lemma~\ref{complexity}.4 and~(\ref{rholis}),
\begin{equation} \label{quiero2}
  \left\| P_{[\ell, u]}\left( x^k - \nabla\left( \|h(x^k)\|^2 + \|g(x^k)_+\|^2\right) \right) -
    x^k \right\| \leq \deltalow.
\end{equation}
 
By (\ref{kendmayor}) and (\ref{kbajo}), we have that $k \leq \kend$,
so, by (\ref{quiero2}),
 \begin{equation} \label{feafea}
\|h(x^k)\|_\infty \leq \delta \mbox{ and } \|g(x^k)_+\|_\infty \leq \delta.
\end{equation}

By (\ref{feafea}), $ g_j(x^k) \leq \delta$ for all $j=1,\dots,p$.
Now, if $g_j(x^k) < -\delta$, we have that $\bar{\mu}^k_j + \rho_k
g_j(x^k) < \bar{\mu}^k_j - \delta \rho_k$, which is smaller than zero
because of~(\ref{rhomuma}), so $\mu_j^{k+1}=0$.

Therefore, the approximate feasibility and complementarity conditions
\begin{equation} \label{feco}
\|h(x^k)\|_\infty \leq \delta,  \; \|g(x^k)_+\| \leq \delta, \mbox{ and }
  \mu^k_j  = 0 \mbox{ if } g_j(x^k) < - \delta
\end{equation}
hold at $x^k$. Moreover, by (\ref{kene}) and
Lemma~\ref{lemcomplexity1}, we have that~(\ref{primaldual1}) also
holds. Therefore, we proved that (\ref{kendmayor}), (\ref{kbajo}),
(\ref{coneps}), (\ref{rhouno}), (\ref{rhocf}), (\ref{rholis}),
(\ref{rhomuma}), and (\ref{kene}) imply (\ref{primaldual1}),
(\ref{primaldual2}), and (\ref{primaldual3}). So, we only need to show
that there exists~$k$ that satisfies
(\ref{kbajo})--(\ref{kene}) or satisfies (\ref{kbajo}),
(\ref{primaldual1}), (\ref{primaldual2}), and (\ref{primaldual3}). In
other words, we must prove that, before completing
\[
N(\deltalow, \varepsilon) + \left[ \frac{\log
    (\delta/\cbig)}{\log(\tau})\right] \times \left[ \frac{\log \left(
    \rho_{\max} / \rho_1 \right)}{\log(\gamma)} \right],
\]   
iterations, we get (\ref{primaldual1}), (\ref{primaldual2}), and
(\ref{primaldual3}) or we get (\ref{kbajo})--(\ref{kene}).

To prove this statement, suppose that, for all $k$ satisfying
(\ref{kbajo}), at least one among the conditions (\ref{primaldual1}),
(\ref{primaldual2}), and (\ref{primaldual3}) does not
hold. Since~(\ref{primaldual1}) necessarily holds if $k \geq
N(\deltalow,\varepsilon)$, this implies that for all~$k$
satisfying~(\ref{kbajo}) and~(\ref{kene}) at least one among the
conditions (\ref{primaldual2}) and (\ref{primaldual3}) does not
hold. By Lemma~\ref{lemcomplexity1}, this implies that for all~$k$
satisfying~(\ref{kbajo}) and~(\ref{kene}),
\[
\max \{\|h(x^k)\|_\infty, \|V_k\|_\infty\} > \delta.   
\]
Then, by (\ref{defcbig}), for $k \geq N(\deltalow, \varepsilon)$, the
existence of more than $\log(\delta/\cbig) / \log(\tau)$ consecutive
iterations $k, k+1, k+2, \dots$ satisfying~(\ref{testfeas})
and~(\ref{kbajo}) is impossible.

Therefore, after the first $N(\deltalow, \varepsilon)$ iterations,
if~$\rho_k$ is increased at iterations $k_1 < k_2$, but not at any
iteration $k \in (k_1, k_2)$, we have that $k_2 - k_1 \leq
\log(\delta/\cbig) / \log(\tau)$. This means that, after the first
$N(\deltalow, \varepsilon)$ iterations, the number of iterations at
which $\rho_k$ is not increased is bounded above by
$\log(\delta/\cbig) / \log(\tau)$ times the number of iterations at
which $\rho_k$ is increased. Now, for obtaining
(\ref{rhouno})--(\ref{rhomuma}), $\log(\rho_{\max}/\rho_1) /
\log(\gamma)$ iterations in which~$\rho_k$ is increased are obviously
sufficient. This completes the desired result.
\end{pro}

\begin{teo} \label{teocomplexity4}
In addition to the hypotheses of Theorem~\ref{teocomplexity3}, assume
that there exist $\cinner > 0$, $v > 0$, and $q > 0$, where $\cinner$
only depends on $\lambda_{\min}$, $\lambda_{\max}$, $\mu_{\max}$,
$\ell$, $u$, and characteristics of the functions $f$, $h$, and $g$,
such that the number of inner iterations, function and derivative
evaluations that are necessary to obtain~(\ref{subprostop}) is bounded
above by $\cinner \rho_k^v \varepsilon_k^{-q}$. Then, the number of
inner iterations, function evaluations, and derivative evaluations
that are necessary to obtain~$k$ such that~(\ref{paradamala}) holds
or~(\ref{primaldual1}), (\ref{primaldual2}) and~(\ref{primaldual3})
hold is bounded above by
\[
\cinner \rho_{\max}^v \varepsilon_{\min,3}^{-q}
\left\{
N(\deltalow, \varepsilon) + \left[ \frac{\log
    (\delta/\cbig)}{\log(\tau})\right] \times \left[ \frac{\log \left(
    \rho_{\max} / \rho_1 \right)}{\log(\gamma)} \right]
\right\},
\]
where $\rho_{\max}$ is given by~(\ref{rhomax}) and 
\begin{equation} \label{epsmin3}
\varepsilon_{\min,3} = \min\{\varepsilon_k \;|\; k \leq  N(\deltalow, \varepsilon) + \left[ \frac{\log
    (\delta/\cbig)}{\log(\tau})\right] \times \left[ \frac{\log \left(
    \rho_{\max} / \rho_1 \right)}{\log(\gamma)} \right].  
\end{equation}
\end{teo}

\begin{pro}
The desired result follows from Theorem~\ref{teocomplexity3} and the
assumptions of this theorem.
\end{pro}

The comparison between Theorems~\ref{teocomplexitygrande2}
and~\ref{teocomplexity4} is interesting. This comparison seems to
indicate that, if we want to be confident that the diagnostic ``$x^k$
is an infeasible stationary point of infeasibility'' is correct, we
must be prepared to pay for that certainty. In fact, the bound
$\rho_{\max}$ on the penalty parameter for the algorithm is defined by
(\ref{rhomax}), which not only grows with $1/\deltalow$, but also
depends on global bounds of the problem $\clips$ and $\cf$. Moreover,
$\varepsilon_k$ also needs to decrease below $\deltalow/4$ because the
decrease of the projected gradient of infeasibility is only guaranteed
by a stronger decrease of the projected gradient of the Augmented
Lagrangian.

\section{Solving the Augmented Lagrangian subproblems} \label{newtonls}

The problem considered in this section is
\begin{equation} \label{theproblem}
\mbox{Minimize } \Phi(x) \mbox{ subject to } x \in \Omega,
\end{equation}
where $\Omega = \{ x \in \R^n \;|\; \ell \leq x \leq u \}$. We assume
that $\Phi$ has continuous first derivatives and that second
derivatives exist almost everywhere. When the Hessian at a point $x$
does not exist we call $\nabla^2 \Phi(x)$ the limit of $\nabla^2
\Phi(x^j)$ for a sequence $x^j$ that converges to $x$.
Problem~(\ref{theproblem}) is of the same type of the problem that is
approximately solved at Step~1 of Algorithm~\ref{al}.1 and we have in
mind the case $\Phi(x) \equiv L_{\rho_k}(x,{\bar \lambda}^k, {\bar
  \mu}^k)$.

For all $I \subseteq \{1, \dots, 2n\}$, we define the 
\textit{open face}
\[
F_I = \{ x \in \Omega \; | \;
x_i = \ell_i \mbox{ if } i \in I, \; x_i = u_i \mbox{ if }
n + i \in I, \; \ell_i < x_i < u_i \mbox{ otherwise} \}.
\]
By definition, $\Omega$ is the union of its open faces and the open
faces are disjoint. This means that every $x \in \Omega$ belongs to
exactly one face $F_I$. The variables $x_i$ such that $\ell_i < x_i <
u_i$ are called \textit{free variables}. For every $x \in \Omega$, we
also define the continuous projected gradient of $\Phi$ given by
\begin{equation} \label{cpg}
g_P(x) = P_{\Omega}(x - \nabla \Phi(x)) - x
\end{equation}
and, if $F_I$ is the open face to which~$x$ belongs, the continuous
projected internal gradient $g_I(x)$ given by
\[
[g_I(x)]_i =
\left\{
\begin{array}{ll}
[g_P(x)]_i, & \mbox{if } x_i \mbox{ is a free variable}, \\
0, & \mbox{otherwise}.
\end{array}
\right.
\]
Note that, sometimes, we write $g_I(x)$ omitting the fact that the
subindex~$I$ refers the face~$F_I$ to which the argument~$x \in
\Omega$ belongs.

The bound-constrained minimization method described in the current
section can be seen as a second-order counterpart of the method
introduced in~\cite{bmgencan}. (See, also, \cite{betra}.)  The
iterates visit the different faces of the box $\Omega$ preserving the
current face while the quotient $\|g_I(x)\|/\|g_P(x)\|$ is big enough
or the new iterate does not hit the boundary. When this quotient
reveals that few progress can be expected from staying in the current
face, the face is abandoned by means of a spectral projected
gradient~\cite{bmr,bmr2,bmr3} iteration. Within each face, iterations
obey a safeguarded Newton scheme with line searches. The employment of
this method is coherent with the conservative point of view of
Algencan. For example, we do not aim to predict the active constraints
at the solution and the inactive bounds have no influence in the
iterations independently of the distance of the current iterate to a
bound. Moreover, we do not try to use second-order information for
leaving the faces. Of course, we do not deny the efficiency of methods
that employ such procedures, but we feel comfortable with the
conservative strategy because the number of algorithmic parameters can
be reduced to a minimum.\\

\noindent
\textbf{Algorithm~\ref{newtonls}.1:} Assume that $x^0 \in \Omega$,
$\Phi_{\target} \in \R$, $r \in (0,1]$, $0 < \tau_1 \leq \tau_2 < 1$,
  $\gamma \in (0,1)$, $\beta \in (0,1)$, $0 < \eta$, $0 <
  \lambda_{\min}^{\SPG} < \lambda_{\max}^{\SPG}$, $0 <
  \sigma_{\chico}$, $0 < \sigma_{\min} \leq \sigma_{\max}$, $0 <
  \underline{h} < \bar h$, $t^{\ext}_{\max} \in \N_{\geq 0}$ are
  given. Initialize $k \leftarrow 0$.

\begin{description}
    
\item[Step 1.] If $\Phi(x^k) \leq \Phi_{\target}$ then stop. Otherwise, if
  $\| g_I(x^k) \|_{\infty} \geq r \| g_P(x^k) \|_{\infty}$ then go to
  Step~2 to perform an \textit{inner-to-the-face iteration using
    Newton with line search} else go to Step~5 to perform a
  \textit{leaving-face iteration using spectral projected gradients
    (SPG)}.

\item[Step 2.] Let $\bar n$ be the number of free variables and let
  $\bar H_k \in \R^{\bar n \times \bar n}$ be the Hessian $\nabla^2 \Phi(x^k)$
  in which rows and columns associated with \textit{non} free
  variables were removed.

\item[Step 2.1.] If $\bar H_k$ is positive definite then set $\sigma
  \leftarrow 0$ and compute $\bar d^k \in \R^{\bar n}$ as the solution
  of $\bar H_k d = - \bar g^k$, where $\bar g^k \in \R^{\bar n}$
  corresponds to $\nabla \Phi(x^k)$ with the components associated with
  the \textit{non} free variables removed, and go to Step~2.3.

\item[Step 2.2.] \textit{Inertia correction}

\item[Step 2.2.1.] If $\sigma^{\ini}$ is undefined then set
  $\sigma^{\ini} \leftarrow
  P_{[\sigma_{\min},\sigma_{\max}]}(\sigma_{\chico} \; h)$, where $h =
  P_{[\underline{h},\bar h]}(\max_{\{i=1,\dots,\bar n\}} \left\{
  \left| [\bar H_k]_{ii} \right| \right\} )$.

\item[Step 2.2.2.] Set $\sigma \leftarrow \sigma^{\ini}$ and while
  $\bar H_k + \sigma I$ is not positive definite do $\sigma \leftarrow
  10 \sigma$.

\item[Step 2.2.3.] Compute $\bar d^k \in \R^{\bar n}$ as a solution of
  $( \bar H_k + \sigma I ) d = - \bar g^k$, where $\bar g^k \in
  \R^{\bar n}$ corresponds to $\nabla \Phi(x^k)$ with the components
  associated with the \textit{non} free variables removed.
  
\item[Step 2.2.4.] While $\| \bar d^k \|_2 > \eta \max\{1, \| \bar
  x^k \|_2 \}$ do $\sigma \leftarrow 10 \sigma$ and redefine $\bar
  d^k$ as the solution of $( \bar H_k + \sigma I ) d = - \bar g^k$. ($\bar
  x^k \in \R^{\bar n}$ corresponds to the free components of $x^k$.)

\item[Step 2.3.] Let $d^k \in \R^n$ be the ``expansion'' of $\bar d^k$
  with $[d^k]_i=0$ if $i$ is a \textit{non}-free-variable index.
  
\item[Step 2.4.] If $\sigma > 0$ then set $\sigma^{\ini} \leftarrow
  P_{[\sigma_{\min},\sigma_{\max}]}(\frac{1}{2} \sigma)$. Otherwise,
  set $\sigma^{\ini} \leftarrow
  P_{[\sigma_{\min},\sigma_{\max}]}(\frac{1}{2} \sigma^{\ini})$.
  
\item[Step 3.] \textit{Line search plus possible projection}

\item[Step 3.1.] Compute $\alpha_{\max}$ as the largest $\alpha > 0$
  such that $x^k + \alpha d^k \in \Omega$. If $\alpha_{\max} \geq 1$
  (i.e. $x^k + d^k \in \Omega$) then skip Step~3.2 below (i.e. go to
  Step~3.3).
  
\item[Step 3.2.] Set $x_{\trial} \leftarrow P_{\Omega}(x^k + d^k)$. If
  $\Phi(x_{\trial}) \leq \Phi_{\target}$ or $\Phi(x_{\trial}) \leq \Phi(x^k)$ then
  set $x^{k+1} = x_{\trial}$ and go to Step~4.

\item[Step 3.3.] Set $t \leftarrow \min\{1, \alpha_{\max}\}$ and
  $x_{\trial} \leftarrow x^k + t d^k$. While $\Phi(x_{\trial}) >
  \Phi_{\target}$ and $\Phi(x_{\trial}) > \Phi(x^k) + t \gamma \langle
  \nabla \Phi(x^k), d^k \rangle$, choose $t_{\new} \in [\tau_1 t, \tau_2
    t]$ and set $t \leftarrow t_{\new}$ and $x_{\trial} \leftarrow x^k
  + t d^k$.

\item[Step 3.4.] Set $x^{k+1} = x_{\trial}$. If $t = \alpha_{\max}$ or
  ( $t=1$ and $\langle \nabla \Phi(x_{\trial}), d^k \rangle > \beta
  \langle \nabla \Phi(x^k), d^k \rangle$ ) then go to Step~4. Otherwise,
  go to Step~6.
  
\item[Step 4.] \textit{Extrapolation}

\item[Step 4.1.] Set $t \leftarrow 1$, $x_{\refe} \leftarrow
  x_{\trial}$, and $x_{\ext} \leftarrow P_{\Omega}(x^k + 2^t (
  x_{\trial} - x^k ))$.

\item[Step 4.2.] While $t \leq t^{\ext}_{\max}$ and
  $\Phi(x_{\ext})<f(x_{\refe})$ and $\Phi(x_{\refe}) > \Phi_{\target}$
  do

  $t \leftarrow t + 1$, $x_{\refe} \leftarrow x_{\ext}$, and $x_{\ext}
  \leftarrow P_{\Omega}(x^k + 2^t ( x_{\trial} - x^k ))$.

\item[Step 4.3.] Reset $x^{k+1} = x_{\refe}$ and go to Step~6.

\item[Step 5.] \textit{Leaving-face SPG iteration}

\item[Step 5.1.] If $k=0$ or $\langle x^k - x^{k-1}, \nabla f(x^k) -
  \nabla f(x^{k-1}) \rangle \leq 0$ then set
  \[
  \lambda_k^{\SPG} = \max \left\{ 1, \|x^k\|_2 / \|g_P(x^k)\|_2 \right\}.
  \]
  Otherwise, set $\lambda_k^{\SPG} = \| x^k - x^{k-1} \|_2^2 / \langle
  x^k - x^{k-1}, \nabla \Phi(x^k) - \nabla \Phi(x^{k-1}) \rangle$.

  In any case, redefine $\lambda_k^{\SPG}$ as $\max\{
  \lambda_{\min}^{\SPG} \min\{ \lambda_k^{\SPG}, \lambda_{\max}^{\SPG}
  \} \}$.

\item[Step 5.2.] Set $t \leftarrow 1$, $x_{\trial} \leftarrow
  P_{\Omega}( x^k - \lambda_k^{\SPG} \nabla \Phi(x^k) )$, and $d^k =
  x_{\trial} - x^k$.

\item[Step 5.3.] While $\Phi(x_{\trial}) > \Phi_{\target}$ and
  $\Phi(x_{\trial}) > \Phi(x^k) + t \gamma \langle \nabla \Phi(x^k), d^k
  \rangle$,

  choose $t_{\new} \in [\tau_1 t, \tau_2 t]$ and set $t \leftarrow
  t_{\new}$ and $x_{\trial} \leftarrow x^k + t d^k$.

\item[Step 5.4.] Set $x^{k+1} =  x_{\trial}$.
  
\item[Step 6.] In the current iteration, the definition of $x^{k+1}$
  implied in the evaluation of $\Phi$ at several points named
  $x_{\trial}$. If, for any of them, we have that $\Phi(x_{\trial}) <
  \Phi(x^{k+1})$ then reset $x^{k+1} = x_{\trial}$. In any case, set
  $k \leftarrow k + 1$ and go to Step~1.

\end{description}

\noindent
\textbf{Remark 1.} At Steps~3.3 and~5.3, interpolation is done with
safeguarded quadratic interpolation. This means that, given
\[
t_{\temp} = - \frac{ \langle \nabla \Phi(x^k), d^k \rangle t^2}
{ 2 ( \Phi(x_{\trial}) - \Phi(x^k) - t \langle \nabla \Phi(x^k), d^k \rangle ) },
\]
if $t_{\temp} \in [\tau_1 t, \tau_2 t]$ then $t_{\new} =
t_{\temp}$. Otherwise, $t_{\new} = \frac{1}{2} t$. This choice
requires $0 < \tau_1 \leq \frac{1}{2} \leq \tau_2 < 1$ instead of
simply $0 < \tau_1 \leq \tau_2 < 1$.

\section{Implementation details and parameters} \label{secimpl}

We implemented Algorithms~\ref{al}.1 and~\ref{newtonls}.1 in
Fortran~90. Implementation is freely available at
\url{http://www.ime.usp.br/~egbirgin/}. Interfaces for solving
user-defined problems coded in Fortran~90 as well as problems from the
CUTEst~\cite{cutest} collection are available. All tests reported
below were conducted on a computer with 3.5 GHz Intel Core i7
processor and 16GB 1600 MHz DDR3 RAM memory, running OS X High Sierra
(version 10.13.6). Codes were compiled by the GFortran compiler of GCC
(version 8.2.0) with the -O3 optimization directive enabled.

\subsection{Augmented Lagrangian method} \label{alparam}

Algorithm~\ref{al}.1 was devised to be applied to a scaled version of
problem~(\ref{theproblem}). Following the Ipopt strategy described
in~\cite[p.46]{ipopt}, in the scaled problem, the objective
function~$f$ is multiplied by
\[
s_f = \max \left\{ 10^{-8}, \frac{100}{\max\{1, \| \nabla f(x^0)
  \|_{\infty}\}} \right\},
\]
each constraint~$h_j$ ($j=1,\dots,m$) is multiplied by
\[
s_{h_j} = \max \left\{ 10^{-8}, \frac{100}{\max\{1, \| \nabla h_j(x^0)
  \|_{\infty}\}} \right\},
\]
and each constraint~$g_j$ ($j=1,\dots,p$) is multiplied by
\[
s_{g_j} = \max \left\{ 10^{-8}, \frac{100}{\max\{1, \| \nabla g_j(x^0)
  \|_{\infty}\}} \right\},
\]
where $x^0 \in \R^n$ is the given initial guess. The scaling is
optional and it is used when the input parameter ``scale'' is set to
``true''. If the parameter is set to ``false'', the original problem, that
corresponds to considering all scaling factors equal to one, is solved.

As stopping criterion, we say that an
iterate $x^k \in [\ell,u]$ with its associated Lagrange multipliers
$\lambda^{k+1}$ and $\mu^{k+1}$ satisfies the main stopping criterion
when
\begin{align}
\max \left\{ \| h(x^k)| \|_{\infty}, \| g(x^k)_+ \|_{\infty} \right\}
&\leq \varepsilon_{\feas}, \label{kkt1}\\
\left\| P_{[\ell,u]}\left(x^k - \left[ s_f \nabla f(x^k) + \sum_{j=1}^m \lambda_j^{k+1} s_{h_j}
\nabla h_j(x^k) + \sum_{j=1}^p \mu_j^{k+1} s_{g_j} \nabla g_j(x^k) \right] \right) - x^k
\right\|_{\infty} &\leq \varepsilon_{\opt}, \label{kkt2}\\
\max_{j=1,\dots,p} \left\{ \min \{ - s_{g_j} g_j(x^k), \mu_j^{k+1} \} \right\} &\leq \varepsilon_{\compl}, \label{kkt3}
\end{align}
where $\varepsilon_{\feas} > 0$, $\varepsilon_{\opt} > 0$, and
$\varepsilon_{\compl} > 0$ are given constants. This means that the
stopping criterion requires \textit{unscaled} feasibility with
tolerance $\varepsilon_{\feas}$ plus \textit{scaled} optimality with
tolerance~$\varepsilon_{\opt}$ and \textit{scaled} complementarity
(measured with the $\min$ function) with tolerance
$\varepsilon_{\compl}$. Note that $x^k \in [\ell,u]$, i.e.\ it
satisfies the bound-constraints with zero tolerance. In addition to
this stopping criterion, Algorithm~\ref{al}.1 also stops if the
penalty parameter~$\rho_k$ reaches the value~$\rho_{\mathrm{big}}$ or
if, in three consecutive iterations, the inner solver that is used at
Step~1 fails at finding a point $x^k \in [\ell,u]$ that
satisfies~(\ref{subprostop}).

In~(\ref{subprostop}) and~(\ref{testfeas}), we consider $\| \cdot \| =
\| \cdot \|_{\infty}$. At Step~2, we consider $\varepsilon_1 =
\sqrt{\varepsilon_{\opt}}$ and $\varepsilon_{k} = \max\{
\varepsilon_{\opt}, 0.1 \varepsilon_{k-1} \}$ for $k > 1$; and, at
Step~3, if $\lambda^{k+1} \in [\lambda_{\min}, \lambda_{\max}]^m$ and
$\mu^{k+1}_i \in [0, \mu_{\max}]^p$ then we set $\bar \lambda^{k+1} =
\lambda^{k+1}$ and $\bar \mu^{k+1} = \mu^{k+1}$. Otherwise, we set
$\bar \lambda^{k+1}=0$ and $\bar \mu^{k+1} = 0$. In the numerical
experiments, we set $\varepsilon_{\feas} = \varepsilon_{\opt} =
\varepsilon_{\compl} = 10^{-8}$, $\rho_{\mathrm{big}} = 10^{20}$,
$\lambda_{\min} = -10^{16}$, $\lambda_{\max} = 10^{16}$, $\mu_{\max} =
10^{16}$, $\gamma = 10$, $\tau = 0.5$, ${\bar \lambda}^1 = 0$, ${\bar
  \mu}^1 = 0$, and
\[
\rho_1 = 10 \max \left\{ 1, \frac{|f(x^0)|}{\max\{ \| h(x^0) \|_2^2 + \| g(x^0)_+ \|_2^2 \}} \right\}.
\]

Two additional strategies complete the implementation of
Algorithm~\ref{al}.1. On the one hand, if Algorithm~\ref{al}.1 fails
at finding a point that satisfies~(\ref{kkt1}), the feasibility
problem~(\ref{fisipro}) is tackled with Algorithm~\ref{newtonls}.1
with the purpose of, at least, finding a feasible point to the
original NLP problem~(\ref{nlp}). On the other hand, at every
iteration~$k$, prior to the subproblem minimization at Step~1,
$(x^{k-1},\lambda^k,\mu^k)$ is used as initial guess to perform ten
iterations of the ``pure'' Newton method (no line search, no inertia
correction) applied to the semismooth KKT system~\cite{mqi,qisun}
associated with problem~(\ref{theproblem}), with dimension $3n+m+p$,
given by
\[
\left(
\begin{array}{c}
\nabla f(x) + \sum_{j=1}^m \lambda_j \nabla h_j(x) + \sum_{j=1}^p \mu_j \nabla g_j(x) - \nu^{\ell} + \nu^u\\
h(x)\\
\min\{ -g(x), \mu \}\\
\min\{ x - \ell, \nu^{\ell} \}\\
\min\{ u - x, \nu^{u} \}\\
\end{array}
\right) =
\left(
\begin{array}{c}
  0\\
  0\\
  0\\
  0\\
  0\\
\end{array}
\right),
\]
where $\nu^{\ell}, \nu^u \in \R^n$ are the Lagrange multipliers
associated with the bound constraints $\ell \leq x$ and $x \leq u$,
respectively. This process is related to the so-called acceleration
process described in~\cite{bmfast} in which a different KKT system was
considered. (See~\cite{bmfast} for details.) The stopping criteria for
the acceleration process are (i) ``the Jacobian of the KKT system has
the 'wrong' inertia'', (ii) ``a maximum of 10 iterations was
reached'', and (iii)
\begin{align}
\max \left\{ \| h(x)| \|_{\infty}, \| g(x)_+ \|_{\infty},
\| (\ell - x)_+ \|_{\infty}, \| (x - u)_+ \|_{\infty} \right\} &\leq \varepsilon_{\feas},\\
\left\| \nabla f(x) + \sum_{j=1}^m \lambda_j \nabla h_j(x) +
\sum_{j=1}^p \mu_j \nabla g_j(x) - \nu^{\ell} + \nu^u \right\|_{\infty} &\leq \varepsilon_{\opt},\\
\max \left\{ \max_{j=1,\dots,p} \left\{[\min \{ - g(x), \mu \}]_j\right\}, \max_{i=1,\dots,n} \left\{[\min \{x-\ell,\nu^{\ell}\}]_i\right\},  \max_{i=1,\dots,n} \left\{[\min \{u-x,\nu^u\}]_i\right\} \right\} &\leq \varepsilon_{\compl}.
\end{align}
Note that criterion (iii) corresponds to satisfying approximate KKT
conditions for the \textit{unscaled} original problem~(\ref{nlp}). On
the other hand, differently from an iterate $x^k \in [\ell,u]$ of
Algorithm~\ref{al}.1 that satisfies
(\ref{kkt1},\ref{kkt2},\ref{kkt3}), a point that satisfies
criterion~(iii) may violate the bound constraints with
tolerance~$\varepsilon_{\feas}$.

If the acceleration process stops satisfying criterion~(i) or~(ii),
everything it was done in the acceleration is discarded and the
iterations of Algorithm~\ref{al}.1 continue. On the other hand, assume
that a point satisfying criterion~(iii) was found by the acceleration
process. If $(x^{k-1},\lambda^k,\mu^k)$
satisfies~(\ref{kkt1},\ref{kkt2},\ref{kkt3}) with half the precision,
i.e.\ with $\varepsilon_{\feas}$, $\varepsilon_{\opt}$, and
$\varepsilon_{\compl}$ substituted by $\varepsilon_{\feas}^{1/2}$,
$\varepsilon_{\opt}^{1/2}$, and $\varepsilon_{\compl}^{1/2}$,
respectively, then we say the acceleration was successful, the point
found by the acceleration process is returned, and the optimization
process stops. On the other hand, if $(x^{k-1},\lambda^k,\mu^k)$ is
far from satisfying~(\ref{kkt1},\ref{kkt2},\ref{kkt3}), we believe the
approximate KKT point the acceleration found may be an undesirable
point. The point is saved for further references, but the optimization
process continues; and the next Augmented Lagrangian subproblem is
tackled by Algorithm~\ref{newtonls}.1 starting from~$x^{k-1}$ and
ignoring the point found by the acceleration process.

\subsection{Bound-constrained minimization method} \label{ubcimpl}

As main stopping criterion of Algorithm~\ref{newtonls}.1, we
considered the condition
\begin{equation} \label{stopcrit}
\| g_P(x^k) \|_{\infty} \leq \varepsilon
\end{equation}
where $g_P(x^k) = P_{[\ell,u]}\left( x^k - \nabla \Phi(x^k) \right) -
x^k$ as defined in~(\ref{cpg}). When an unconstrained or
bound-constrained problem is being solved, in~(\ref{stopcrit}) and in
the alternative stopping criteria described below, we use $\varepsilon
= \varepsilon_{\opt} = 10^{-8}$. When the problem being tackled by
Algorithm~\ref{newtonls}.1 is a subproblem of Algorithm~\ref{al}.1,
the value of $\varepsilon$ in~(\ref{stopcrit}) and in the alternative
stopping criteria described below is the one described in
Section~\ref{alparam} (that we cannot mention here since we are using
$k$ to denote iterations of both Algorithms~\ref{al}.1
and~\ref{newtonls}.1). In addition, Algorithm~\ref{newtonls}.1 may
also stop at iteration~$k$ by any of the following alternative
stopping criteria: \textbf{(a)}
$\|g_P(x^{k-\ell})\|_{\infty}<\sqrt{\varepsilon}$ for all $0 \leq \ell
< 100$; \textbf{(b)} $\|g_P(x^{k-\ell})\|_{\infty}<\varepsilon^{1/4}$
for all $0 \leq \ell < 5{,}000$; \textbf{(c)}
$\|g_P(x^{k-\ell})\|_{\infty}<\varepsilon^{1/8}$ for all $0 \leq \ell
< 10{,}000$; \textbf{(d)} $\Phi(x^k) \leq \Phi_{\target}$;
\textbf{(e)} $k \geq k_{\max} = 50{,}000$; and \textbf{(f)}
$k_{\best}$ is the smallest index such that $\Phi(x^{k_{\best}}) =
\min \{ \Phi(x^0), \Phi(x^1), \dots, \Phi(x^k) \}$ and $k - k_{\best}
> 3$, i.e. the best functional value so far obtained is not updated in
three consecutive iterations.

In Algorithm~\ref{newtonls}.1, although the theory allows us a wide
range of possibilities, in practice we consider $H_k = \nabla^2
\Phi(x^k)$ for all~$k$.  The linear systems at Step~2.1 and~2.2.2 are
solved with subroutine MA57 from HSL~\cite{hsl} (using all its default
parameters). In the experiments, we set $\Phi_{\target}=-10^{12}$,
$r=0.1$, $\tau_1=0.1$, $\tau_2=0.9$, $\gamma=10^{-4}$, $\beta=0.5$,
$\eta=10^4$, $\lambda_{\min}^{\SPG} = 10^{-16}$,
$\lambda_{\max}^{\SPG} = 10^{16}$, $\sigma_{\chico} = 10^{-8}$,
$\sigma_{\min} = 10^{-8}$, $\sigma_{\max} = 10^{16}$, $\underline{h} =
10^{-8}$, $\bar h = 10^{8}$, and $t^{\ext}_{\max} = 20$.

When Algorithm~\ref{newtonls}.1 is used to solve a subproblem of
Algorithm~\ref{al}.1, we have that $\nabla^2 \Phi(x) = \nabla^2
L_{\rho_k}(x,{\bar \lambda}_k,{\bar \mu}_k)$, i.e. $\nabla^2 \Phi(x)$
is the Hessian of the augmented Lagrangian associated with the scaled
version of problem~(\ref{theproblem}) given by
\begin{equation} \label{lahess}
\resizebox{\textwidth}{!}{$
s_f \nabla^2 f(x) +
\sum_{j=1}^m \left\{ {\bar \lambda}_j^k s_{h_j} \nabla^2 h_j(x) + \rho_k s_{h_j}^2 \nabla h_j(x) \nabla h_j(x)^T \right \} +
\sum_{j \in I_k} \left\{ {\bar \mu}_j^k s_{g_j} \nabla^2 g_j(x) + \rho_k s_{g_j}^2 \nabla g_j(x) \nabla g_j(x)^T \right\},
$}
\end{equation}
where $I_k = I_{\rho_k}(x^k,{\bar \mu}^k) = \{ j = 1,\dots,p \;|\;
{\bar \mu}^k + \rho_k s_{g_j} g_j(x^k) > 0 \}$. A relevant issue from
the practical point of view is that, despite the sparsity of the
Hessian of the Lagrangian and the sparsity of the Jacobian of the
constraints, this matrix may be dense. Thus, factorizing, or even
building it, may be prohibitive. As an alternative, instead of
building and factorizing the Hessian above, it can be solved an
augmented linear system with the coefficients' matrix given by
\begin{equation} \label{lahess2}
\left(
\begin{array}{c|c}
s_f \nabla^2 f(x) +
\sum_{j=1}^m \left\{ {\bar \lambda}_j^k s_{h_j} \nabla^2 h_j(x) \right \} +
\sum_{j \in I_k} \left\{ {\bar \mu}_j^k s_{g_j} \nabla^2 g_j(x) \right\} & J(x)^T \\[2mm]
\hline
\phantom{\displaystyle \sum} J(x) & - \frac{1}{\rho_k} I
\end{array}
\right),
\end{equation}
where $J(x)$ is a matrix whose columns are $\nabla h_1(x), \dots,
\nabla h_m(x)$ plus the gradients $\nabla g_j(x)$ such that $j \in
I_k$. This matrix preserves the sparsity of the Hessian of the
Lagrangian and of the Jacobian of the constraints. The implementation
of Algorithms~\ref{newtonls}.1 dynamically selects one of the two
aproaches.

Another relevant fact from the practical point of view, related to
matrices~(\ref{lahess}) and~(\ref{lahess2}), is that the current tools
available in CUTEst compute the full Jacobian of the constraints and
$\sum_{j=1}^p {\bar \mu}_j^k s_{g_j} \nabla^2 g_j(x)$ with ${\bar
  \mu_j}^k=0$ if $j \not\in I_k$ instead of $J(x)$ and $\sum_{j \in
  I_k} {\bar \mu}_j^k s_{g_j} \nabla^2 g_j(x)$, respectively. On the
one hand, this feature preserves the Jacobian's and the
Hessian-of-the-Lagrangian's sparsity structures independently
of~${\bar \mu^k}$ and~$x$, as required by some solvers. On the other
hand, it impairs Algorithm~\ref{al}.1, when applied to problems from the
CUTEst collection, of fully exploiting the potential advantage of
dealing with inequality constraints without adding slack variables. In
summary, Algorithm~\ref{al}.1--\ref{newtonls}.1 is prepared to deal
with matrices with different sparsity structures at every iteration
and, for that reason, it performs the analysis step of the
factorization at every iteration. This is the price to pay for
exploiting inequality constraints without adding slack
variables. However, the CUTEst subroutines are not prepared to exploit
this feature and Algorithm~\ref{al}.1--\ref{newtonls}.1, when solving
problems from the CUTEst collection, pays the price without enjoying
the advantages. Of course, this CUTEst inconvenient influences
negatively the comparison of Algencan with other solvers if the CPU
time is used as a performance measure.

\section{Numerical experiments} \label{secnumexp}

In this section, we aim to evaluate the performance of
Algorithm~\ref{al}.1--\ref{newtonls}.1 (referred as Algencan from now
on) for solving unconstrained, bound-constrained, feasibility, and
nonlinear programming problems. The performance of Ipopt~\cite{ipopt}
(version 3.12.12) is also exhibited. Both methods were run in the same
computational environment, compiled with the same BLAS routines, and
also using the same subroutine MA57 from HSL for solving the linear
systems. All Ipopt default parameters were used\footnote{Option
  'honor\_original\_bounds no', that does not affect Ipopt's
  optimization process, was used. Ipopt might relax the bounds during
  the optimization beyond its initial \textit{relative} relaxation
  factor whose default value is $10^{-8}$. Option
  'honor\_original\_bounds no' simply avoids the final iterate to be
  projected back onto the box defined by the bound constraints. So,
  the actual absolute violation of the bound constraints at the final
  iterate can be measured.}. A CPU time limit of 10 minutes per
problem was imposed. In the numerical experiments, we considered all
$1{,}258$ problems from the CUTEst collection~\cite{cutest} with their
default dimensions. In the collection, there are 217 unconstrained
problems, 144 bound-constrained problems, 157 feasibility problems,
and 740 nonlinear programming problems. A hint on the number of
variables in each family is given in Table~\ref{tab0}.

\begin{table}[h!]
\begin{center}
\begin{tabular}{cccccc}
\hline
\multirow{2}{*}{Problem type} & \multirow{2}{*}{\# of problems} & \multirow{2}{*}{$n_{\max}$} &
\multicolumn{3}{c}{\# of problems with $n \geq \omega n_{\max}$} \\
\cline{4-6}
& & & $\omega=0.1$ & $\omega=0.01$ & $\omega=0.001$ \\
\hline
unconstrained     & 217 & 100{,}000 & 15 &  87 &  97 \\
bound-constrained & 144 & 149{,}624 &  5 &  60 &  72 \\
feasibility       & 156 & 123{,}200 &  5 &  40 &  55 \\
NLP               & 740 & 250{,}997 & 67 & 263 & 379 \\
\hline
\end{tabular}
\end{center}
\caption{Distribution of the number of variables~$n$ in the CUTEst
  collection test problems.}
\label{tab0}
\end{table}

Large tables with a detailed description of the output of each method
in the $1{,}258$ problems can be found in
\url{http://www.ime.usp.br/~egbirgin/}. A brief overview follows. Note
that, since the methods differ in the stopping criteria, arbitrary
decisions will be made. A point in common is that both methods seek
satisfying the (sup-norm of the) violation of the unscaled equality
and inequality constraints with precision~$\varepsilon_{\feas} =
10^{-8}$. However, as described in~\cite[\S3.5]{ipopt}, Ipopt
considers a \textit{relative} initial relaxation of the bound
constraints (whose default value is $10^{-8}$); and it may apply
repeated additional relaxations during the optimization
process. Table~\ref{tabfeas} shows the number of problems in which
each method found a point satisfying
\begin{equation} \label{feas1}
\max\{ \| h(x) \|_{\infty}, \| [g(x)]_+ \|_{\infty} \} \leq \varepsilon_{\feas}
\end{equation}
plus
\begin{equation} \label{feas2}
\max\{ \| (\ell - x)_+ \|_{\infty}, \| (x - u)_+ \|_{\infty} \} \leq \bar \varepsilon_{\feas}
\end{equation}
with $\varepsilon_{\feas} = 10^{-8}$ and $\bar \varepsilon_{\feas} \in
\{ 0.1, 10^{-2}, \dots, 10^{-16}, 0 \}$. Figures in the table show
that, in most cases, Algencan satisfies the bound constraints
with zero tolerance and that the violation of the bound constraints
 rarely exceeds the tolerance~$10^{-8}$. This is an expected result,
since the method satisfies these requirements by definition. Regarding
Ipopt, the table shows in which way the amount of problems in
which~(\ref{feas2}) holds varies as a function of the tolerance~$\bar
\varepsilon_{\feas}$.

\begin{table}[ht!]
\begin{center}
\resizebox{\textwidth}{!}{
\begin{tabular}{lrrrrrrrrrrrrrrrrr}
\cline{2-18}
 & \multicolumn{17}{c}{$\bar \varepsilon_{\feas}$}\\
\hline
 & $0.1$ & $10^{-2}$ & $10^{-3}$ & $10^{-4}$ & $10^{-5}$ & $10^{-6}$ & $10^{-7}$ & $10^{-8}$ 
 & $10^{-9}$ & $10^{-10}$ & $10^{-11}$ & $10^{-12}$ & $10^{-13}$ & $10^{-14}$ & $10^{-15}$ & $10^{-16}$ & $0$ \\  
\hline
Algencan &
1{,}132 & 1{,}132 & 1{,}131 & 1{,}131 & 1{,}131 & 1{,}130 & 1{,}130 & 1{,}130 & 1{,}121 & 1{,}115 & 1{,}112 & 1{,}105 & 1{,}092 & 1{,}082 & 1{,}077 & 1{,}069 & 1{,}058\\
Ipopt &
1{,}073 & 1{,}072 & 1{,}070 & 1{,}068 & 1{,}056 & 1{,}044 & 1{,}016 &  970 &  794 &  793 &  793 &  793 &  793 &  792 &  792 &  792 &  791\\
\hline
\end{tabular}}
\end{center}
\caption{Number of problems in which a point
  satisfying~(\ref{feas1},\ref{feas2}) was found by Algencan and Ipopt
  with $\varepsilon_{\feas} = 10^{-8}$ and $\bar \varepsilon_{\feas}
  \in \{ 0.1, 10^{-2}, \dots, 10^{-16}, 0\}$.}
\label{tabfeas}
\end{table}

If the violation of the bound constraints is disregarded,
Table~\ref{tabfeas} shows that Algencan found points
satisfying~(\ref{feas1},\ref{feas2}) with $\varepsilon_{\feas} =
10^{-8}$ and $\bar \varepsilon_{\feas}=0.1$ in~$1{,}132$ problems;
while Ipopt found the same in~$1{,}073$. There are in the CUTEst
collection 85 problems (62 feasibility problems and~23 nonlinear
programming problems) in which the number of equality constraints is
larger than the number of variables. Ipopt does not apply to these
problems and, thus, of course, it does not find a point
satisfying~(\ref{feas1},\ref{feas2}). Algencan \textit{did} find a
point satisfying~(\ref{feas1},\ref{feas2}) in~28 out of the~85
problems to which Ipopt does not apply; and this explains half of the
difference between the methods. In any case, it can be said that, over
a universe of $1{,}258$ problems, both methods found ``feasible
points'' in a large fraction of the problems; recalling that the
collection contains infeasible problems.

We now consider the set of~757 problems in which both methods found a
point satisfying~(\ref{feas1}) with $\varepsilon_{\feas} = 10^{-8}$
and~(\ref{feas2}) with $\bar \varepsilon_{\feas} = 0$. For a given
problem, let $f_1$ be the value of the objective function at the point
found by Algencan; let $f_2$ be the value of the objective function at
the point found by Ipopt; and let $f^{\min} = \min\{ f_1, f_2
\}$. Table~\ref{tabbest} shows in how many problems it holds
\begin{equation} \label{eqbest}
f_i \leq f^{\min} + f_{\tol} \max\{ 1, | f^{\min} | \} \mbox{ for } i = 1, 2
\end{equation}
and $f_{\tol} \in \{ 0.1, 10^{-2}, \dots, 10^{-8}, 0 \}$.

\begin{table}[ht!]
\begin{center}
\begin{tabular}{lrrrrrrrrr}
\cline{2-10}
 & \multicolumn{9}{c}{$f_{\tol}$}\\
\hline
 & $0.1$ & $10^{-2}$ & $10^{-3}$ & $10^{-4}$ & $10^{-5}$ & $10^{-6}$ & $10^{-7}$ & $10^{-8}$ & 0 \\
\hline
Algencan & 722 & 715 & 706 & 694 & 691 & 678 & 675 & 663 & 498 \\
Ipopt    & 723 & 708 & 699 & 694 & 683 & 653 & 623 & 592 & 383 \\
\hline
\end{tabular}
\end{center}
\caption{Number of problems in which a point satisfying~(\ref{feas1})
  with $\varepsilon_{\feas} = 10^{-8}$, (\ref{feas2}) with $\bar
  \varepsilon_{\feas} = 0$, and~(\ref{eqbest}) with $f_{\tol} \in \{
  0.1, 10^{-2}, \dots, 10^{-8}, 0 \}$ was found by Algencan and
  Ipopt.}
\label{tabbest}
\end{table}

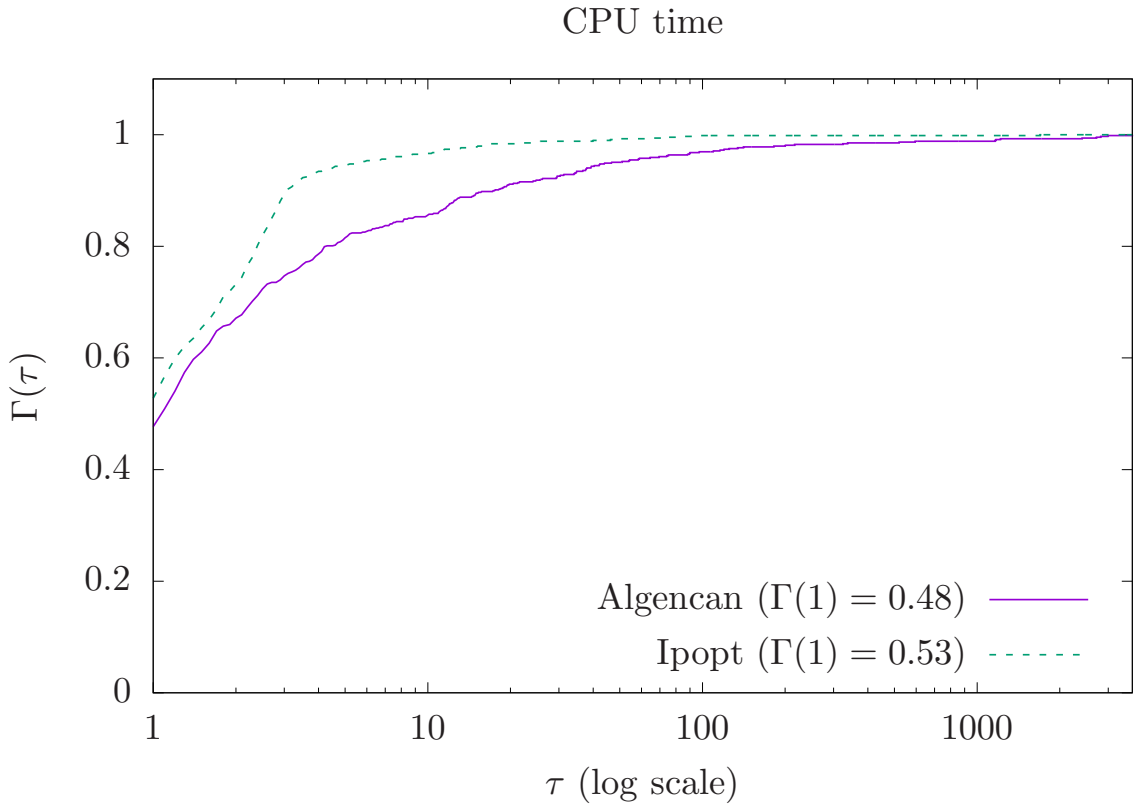
\begin{figure}[h!]
\begin{center}
\resizebox{1.0\textwidth}{!}{\input{bmcomperfigpptime.tex}}
\end{center}
\caption{Performance profiles comparing the CPU time spent by Algencan
  and Ipopt in the~688 problems in which both methods found a point
  that satisfies~(\ref{feas1}) with $\varepsilon_{\feas} = 10^{-8}$,
  (\ref{feas2}) with $\bar \varepsilon_{\feas} = 0$,
  and~(\ref{eqbest}) with $f_{\tol} = 0.1$.}
\label{figpp}
\end{figure}

Finally, we consider the set of~688 problems in which both, Algencan
and Ipopt, found a point that satisfies~(\ref{feas1}) with
$\varepsilon_{\feas} = 10^{-8}$, (\ref{feas2}) with $\bar
\varepsilon_{\feas} = 0$, and~(\ref{eqbest}) with $f_{\tol} =
0.1$. For this set of problems, Figure~\ref{figpp} shows the
performance profile~\cite{pp} that considers, as performance measure,
the CPU time spent by each method. In the figure, for $i \in M \equiv
\{ \mathrm{Algencan}, \mathrm{Ipopt} \}$,
\[
\Gamma_i(\tau) = \frac{\#\left\{ j \in \{1,\dots,q\} \; | \; t_{ij} \leq \tau \min_{s \in M} \{ t_{sj} \} \right\}}{q},
\]
where $\#{\cal S}$ denotes the cardinality of set~${\cal S}$, $q=688$
is the number of considered problems, and $t_{ij}$ is the performance
measure (CPU time) of method~$i$ applied to problem~$j$. Thus,
$\Gamma_{\mathrm{Algencan}}(1) = 0.48$ and $\Gamma_{\mathrm{Ipopt}}(1)
= 0.53$ says that Algencan was faster than Ipopt in 48\% of the
problems and Ipopt was faster then Algencan in 53\% of the problems.
Complementing the performance profile, we can report that there are~9
problems in which both methods spent at least a second of CPU time and
one of the methods is at least ten times faster than the other. Among
these~9 problems, Ipopt is faster in~5 and Algencan is faster in the
other~4.

\section{Conclusions} \label{secconcl}

In this work, a version of the (safeguarded) Augmented Lagrangian
algorithm Algencan \cite{abmstango,bmbook} that possesses iteration
and evaluation complexity was described, implemented, and
evaluated. Moreover, the convergence theory of Algencan was
complemented with new complexity results. The way in which an
Augmented Lagrangian method was able to inherit the complexity
properties from a method for bound-constrained minimization is a nice
example of the advantages of the modularity feature that Augmented
Lagrangian methods usually possess.

As a byproduct of this development, a new version of Algencan that
uses a Newtonian method with line search to solve the subproblems was
developed from scratch. Moreover, the acceleration process described
in~\cite{bmfast} was revisited. In particular, the KKT system with
complementarity modelled with the product between constraints and
multipliers was replaced with the KKT system that models the
complementarity constraints with the semismooth $\min$ function.
  
We provided a fully reproducible comparison with Ipopt, which is,
probably, de most effective and best known free software for
constrained optimization. The main feature we want to stress is that
there exist a significative number of problems that Algencan solves
satisfactorily whereas Ipopt does not, and vice versa. This is not
surprising because the way in which Augmented Lagrangians and Interior
Point Newtonian methods handle problems are qualitatively
different. Constrained Optimization is an extremely heterogeneous
family. Therefore, we believe that what justifies the existence of new
algorithms or the survival of traditional ones is not their capacity
of solving a large number of problems using slightly smaller computer
time than ``competitors'', but the potentiality of solving some
problems that other algorithms fail to solve. Engineers and
practitioners should not care about the choice between algorithm~A
or~B according to subtle efficiency criteria. The best strategy is to
contemplate both, using one or the other according to their behavior
on the family of problems that they need to solve in practice. As in
many aspects of life, competition should give place to cooperation. \\

\noindent
\textbf{Acknowledgements.} The authors are indebted to Iain Duff, Nick
Gould, Dominique Orban, and Tyrone Rees for their help in issues
related to the usage of MA57 from HSL and the CUTEst collection.

\end{document}

%% file: bmcomperfigpptime.tex
\begingroup
  \makeatletter
  \providecommand\color[2][]{%
    \GenericError{(gnuplot) \space\space\space\@spaces}{%
      Package color not loaded in conjunction with
      terminal option `colourtext'%
    }{See the gnuplot documentation for explanation.%
    }{Either use 'blacktext' in gnuplot or load the package
      color.sty in LaTeX.}%
    \renewcommand\color[2][]{}%
  }%
  \providecommand\includegraphics[2][]{%
    \GenericError{(gnuplot) \space\space\space\@spaces}{%
      Package graphicx or graphics not loaded%
    }{See the gnuplot documentation for explanation.%
    }{The gnuplot epslatex terminal needs graphicx.sty or graphics.sty.}%
    \renewcommand\includegraphics[2][]{}%
  }%
  \providecommand\rotatebox[2]{#2}%
  \@ifundefined{ifGPcolor}{%
    \newif\ifGPcolor
    \GPcolorfalse
  }{}%
  \@ifundefined{ifGPblacktext}{%
    \newif\ifGPblacktext
    \GPblacktexttrue
  }{}%
  \let\gplgaddtomacro\g@addto@macro
  \gdef\gplbacktext{}%
  \gdef\gplfronttext{}%
  \makeatother
  \ifGPblacktext
    \def\colorrgb#1{}%
    \def\colorgray#1{}%
  \else
    \ifGPcolor
      \def\colorrgb#1{\color[rgb]{#1}}%
      \def\colorgray#1{\color[gray]{#1}}%
      \expandafter\def\csname LTw\endcsname{\color{white}}%
      \expandafter\def\csname LTb\endcsname{\color{black}}%
      \expandafter\def\csname LTa\endcsname{\color{black}}%
      \expandafter\def\csname LT0\endcsname{\color[rgb]{1,0,0}}%
      \expandafter\def\csname LT1\endcsname{\color[rgb]{0,1,0}}%
      \expandafter\def\csname LT2\endcsname{\color[rgb]{0,0,1}}%
      \expandafter\def\csname LT3\endcsname{\color[rgb]{1,0,1}}%
      \expandafter\def\csname LT4\endcsname{\color[rgb]{0,1,1}}%
      \expandafter\def\csname LT5\endcsname{\color[rgb]{1,1,0}}%
      \expandafter\def\csname LT6\endcsname{\color[rgb]{0,0,0}}%
      \expandafter\def\csname LT7\endcsname{\color[rgb]{1,0.3,0}}%
      \expandafter\def\csname LT8\endcsname{\color[rgb]{0.5,0.5,0.5}}%
    \else
      \def\colorrgb#1{\color{black}}%
      \def\colorgray#1{\color[gray]{#1}}%
      \expandafter\def\csname LTw\endcsname{\color{white}}%
      \expandafter\def\csname LTb\endcsname{\color{black}}%
      \expandafter\def\csname LTa\endcsname{\color{black}}%
      \expandafter\def\csname LT0\endcsname{\color{black}}%
      \expandafter\def\csname LT1\endcsname{\color{black}}%
      \expandafter\def\csname LT2\endcsname{\color{black}}%
      \expandafter\def\csname LT3\endcsname{\color{black}}%
      \expandafter\def\csname LT4\endcsname{\color{black}}%
      \expandafter\def\csname LT5\endcsname{\color{black}}%
      \expandafter\def\csname LT6\endcsname{\color{black}}%
      \expandafter\def\csname LT7\endcsname{\color{black}}%
      \expandafter\def\csname LT8\endcsname{\color{black}}%
    \fi
  \fi
    \setlength{\unitlength}{0.0500bp}%
    \ifx\gptboxheight\undefined%
      \newlength{\gptboxheight}%
      \newlength{\gptboxwidth}%
      \newsavebox{\gptboxtext}%
    \fi%
    \setlength{\fboxrule}{0.5pt}%
    \setlength{\fboxsep}{1pt}%
\begin{picture}(7200.00,5040.00)%
    \gplgaddtomacro\gplbacktext{%
      \csname LTb\endcsname
      \put(814,704){\makebox(0,0)[r]{\strut{}$0$}}%
      \put(814,1372){\makebox(0,0)[r]{\strut{}$0.2$}}%
      \put(814,2040){\makebox(0,0)[r]{\strut{}$0.4$}}%
      \put(814,2709){\makebox(0,0)[r]{\strut{}$0.6$}}%
      \put(814,3377){\makebox(0,0)[r]{\strut{}$0.8$}}%
      \put(814,4045){\makebox(0,0)[r]{\strut{}$1$}}%
      \put(946,484){\makebox(0,0){\strut{}$1$}}%
      \put(2589,484){\makebox(0,0){\strut{}$10$}}%
      \put(4232,484){\makebox(0,0){\strut{}$100$}}%
      \put(5875,484){\makebox(0,0){\strut{}$1000$}}%
    }%
    \gplgaddtomacro\gplfronttext{%
      \csname LTb\endcsname
      \put(198,2541){\rotatebox{-270}{\makebox(0,0){\strut{}$\Gamma(\tau)$}}}%
      \put(3874,154){\makebox(0,0){\strut{}$\tau$ ($\log$ scale)}}%
      \put(3874,4709){\makebox(0,0){\strut{}CPU time}}%
      \csname LTb\endcsname
      \put(5816,1262){\makebox(0,0)[r]{\strut{}Algencan ($\Gamma(1) = 0.48$)}}%
      \csname LTb\endcsname
      \put(5816,932){\makebox(0,0)[r]{\strut{}Ipopt    ($\Gamma(1) = 0.53$)}}%
    }%
    \gplbacktext
    \put(0,0){\includegraphics{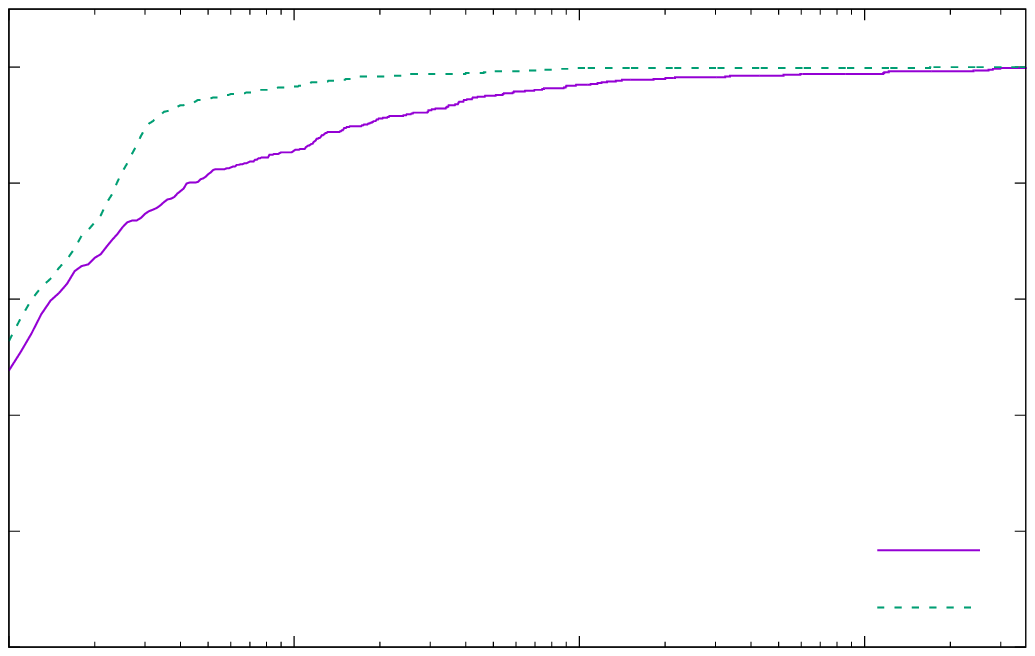}}%
    \gplfronttext
  \end{picture}%
\endgroup